\documentclass[graybox]{svmult}

\usepackage{mathptmx}       
\usepackage{helvet}         
\usepackage{courier}        
\usepackage{type1cm}        
\usepackage{makeidx}         
\usepackage{graphicx}        
\usepackage{multicol}        
\usepackage[bottom]{footmisc}

\usepackage{amsmath}
\usepackage{amsfonts}
\usepackage{amssymb}
\usepackage{color}
\usepackage{epsfig}

\def\mm#1{\begin{align}#1\end{align}}
\def\mmn#1{\begin{align*}#1\end{align*}}

\def\no{{m-1}}
\def\nn{{m}}
\def\dtn{\Delta t_\nn}

\def\R{\mathbb{R}}
\def\Rd1{\Omega\times\R_+^0}

\def\ff{\mathbf{f}}
\def\fn{{\ff_\nu}}
\def\FF{\mathcal{F}}
\def\FFn{{\FF_\nu^\nn}}

\def\At{{\tilde A}}

\def\Pt{{\tilde P}}
\def\inn{{i,\nn}}

\def\O{{\Omega}}
\def\OT{{\Omega_T}}
\def\Ob{{\overline\Omega}}
\def\OTb{{\overline\OT}}
\def\GT{{\Gamma_T}}

\def\Oh{{\Omega_h}}
\def\OTh{{\Omega_{T,h}}}
\def\Oin{{V_i^\nn}}
\def\Ojn{{V_j^\nn}}

\def\Eint{{\mathcal E_{T,h}^{int}}}
\def\Eext{{\mathcal E_{T,h}^{ext}}}
\def\Gin{{\partial V_i^\nn}}
\def\Gij{{\Gamma_{ij}}}
\def\nij{{n_{ij}}}
\def\Gijn{{\Gamma_{ij}^\nn}}

\def\In{{I_\nn}}
\def\Vi{{V_i}}
\def\Vin{{V_i^\nn}}
\def\UU{\mathcal U}
\def\VV{\mathcal V}
\def\VVh{{{\mathcal V}_h}}
\def\myVVh{{{\mathcal V}_h}}
\def\VVh0{{\mathcal V_h^0}}

\def\eJ{{\varepsilon_J}}
\def\eO{{\varepsilon_\Omega}}
\def\egm{{\varepsilon_\Gamma^-}}
\def\egp{{\varepsilon_\Gamma^+}}

\def\vp{\varphi}

\def\ebkn{\bar\eta_k^\nn}

\def\tol{T\!ol}
\def\cfl{C\!F\!L}
\def\myN{{\mathcal N}}

\makeindex             

\begin{document}

\title*{Timestep control for weakly instationary flows}
\author{Christina Steiner and Sebastian Noelle}
\institute{Christina Steiner \at Institut f\"ur Geometrie und Praktische Mathematik,
  RWTH Aachen University, Templergraben 55, D-52056 Aachen, Germany.
        \email{steiner@igpm.rwth-aachen.de}
\and Sebastian Noelle \at Institut f\"ur Geometrie und Praktische Mathematik,
  RWTH Aachen University, Templergraben 55, D-52056 Aachen, Germany.
        \email{noelle@igpm.rwth-aachen.de}}
%
%

\maketitle

\abstract{We report on recent work on adaptive timestep control for weakly instationary gas flows 
\cite{Steiner_2008_diss, SteinerNoelle, SteinerMuellerNoelle} carried out within SFB 401, TPA3.
The method which we implement and extend is a space-time
splitting of adjoint error representations for target functionals due to S\"uli \cite{Sueli} and Hartmann \cite{Hartmann1998}. In this paper, we first review the method for scalar, 1D, conservation laws.
We design a test problem for weakly instationary solutions and show numerical experiments which clearly show
the possible benefits of the method. Then we extend the approach to the 2D Euler equations of gas dynamics.
New ingredients are (i) a conservative formulation of the adjoint problem which makes its solution robust
and efficient, (ii) the derivation of boundary conditions for this new formulation of the  adjoint problem and (iii) the coupling of the adaptive time-stepping with the multiscale spatial adaptation due to M\"uller \cite{Mueller:02,BramkampLambyMueller2004}, also developed within SFB 401.
The combined space-time adaptive method provides an efficient choice of timesteps for implicit
computations of weakly instationary flows. The timestep will be very large in
regions of stationary flow, and becomes small when a perturbation enters the
flow field. The efficiency of the Euler solver is investigated by means of an unsteady
inviscid 2D flow over a bump.}

\section{Introduction}

For the aeroelastic problems studied in the SFB~401 project both the stationary and the instationary case are of interest. In either case, {\em adaptive spatial grids} are well-established and help to reduce computational time and storage.  There has been a tremendous amount of research designing, analyzing and
implementing codes which are adaptive in space, see e.g.
\cite{clawpack,Mueller:02,KroenerOhlberger2000, Ohlberger2001} and references
therein.

Adaptive time-marching towards the stationary solution is an essential ingredient of most CFD codes. Here the timestep-sizes are adapted to the local CFL numbers. So far, it is not well-understood if adaptive timesteps might be useful for instationary flows. In the present paper we report on the work of SFB~401's project TPA3b, where we explored {\em adaptive explicit/implicit timesteps} for {\em weakly instationary flows}.

While time accuracy is still needed to study
phenomena like aero-elastic interactions, large timesteps may be possible when
the perturbations have passed.  For explicit calculations of instationary
solutions to hyperbolic conservation laws, the timestep is dictated by the
$\cfl$ condition due to Courant, Friedrichs and Lewy
\cite{CourantFriedrichsLewy1928}, which requires that the numerical speed of
propagation should be at least as large as the physical one. For implicit
schemes, the $\cfl$ condition does not provide a restriction, since the
numerical speed of propagation is infinite.  Depending on the equations and the
scheme, restrictions may come in via the stiffness of the resulting nonlinear
problem. These restrictions are usually not as strict as in the explicit case,
where the $\cfl$ number should be below unity. For implicit calculations, $\cfl$
numbers of much larger than 1 may well be possible. Therefore, it is a serious
question how large the timestep, i.e. the $\cfl$ number, should be chosen. 

Some time adaptation stategies developed by other authors are those of Ferm and L\"otstedt
\cite{Loetstedt} based on  timestep control strategies for ODEs, and extended to fully adaptive multiresolution finite volume schemes, see \cite{Domingues2007,Domingues2008}. Alternatively  Kr\"oner and Ohlberger
\cite{KroenerOhlberger2000, Ohlberger2001} based their space-time adaptivity 
upon Kuznetsov-type a-posteriori $L^1$ error-estimates for scalar conservation
laws.

In this paper we will use a space-time-split adjoint error representation to
control the timestep adaptation. For this purpose, let us briefly summarize the
space-time splitting of the adjoint error representation, see
\cite{Johnson4,BeckerRannacher1996,BeckerRannacher2000,Sueli,Hartmann1998}
for details. The error representation expresses the error in a target functional
as a scalar product of the finite element residual with the dual solution. This error
representation is decomposed into separate spatial and temporal components. The
spatial part will decrease under refinement of the spatial grid, and the
temporal part under refinement of the timestep. Technically, this decomposition
is achieved by inserting an additional projection. Usually, in the error
representation, one subtracts from the dual solution its projection onto
space-time polynomials. Now, we also insert the projection of the dual solution
onto polynomials in time having values which are $H^1$ functions with respect to
space.

This splitting can be used to develop a strategy for a local choice of
timestep. In contrast to the results reported in \cite{Sueli,Hartmann1998} for 
scalar conservation laws we now investigate weakly instationary
solution to the 2D Euler equations. The timestep will be very large in regions of stationary flow, and becomes small when a perturbation enters the flow field. 

Besides applying well-established adjoint techniques to a new test problem, we
further develop a new technique (first proposed by the authors in
\cite{SteinerNoelle}) which simplifies and accelerates the computation of the
dual problem. Due to Galerkin orthogonality, the dual solution $\vp$ does not
enter the error representation as such. Instead, the relevant term is the
difference between the dual solution and its projection to the finite element
space, $\vp-\vp_h$. In \cite{SteinerNoelle} we showed that it is therefore
sufficient to compute the spatial gradient of the dual solution, $w=\nabla \vp$.
This gradient satisfies a conservation law instead of a transport equation, and
it can therefore be computed with the same conservative algorithm as the forward
problem \cite{SteinerNoelle}. The great advantage is that the conservative
backward algorithm can handle possible discontinuities in the coefficients
robustly.

A key step is to formulate boundary conditions for the gradient  $w = \nabla
\vp$ instead of $v$. Generally the boundary conditions for the dual problem
come from the weighting functions of the target functional, e.g. lift or drag.
To formulate boundary conditions for $w$ which are compatible with the target
functional, one has to lift the well-established techniques of characteristic
decompositions from the dual solution to its gradient. We will present details
on that in Section~\ref{sect:conservative_dual_problem}.

Starting with a very coarse, but adaptive spatial mesh and $\cfl$ below unity, we
establish timesteps which are well adapted to the physical problem at hand. The
scheme detects stationary time regions, where it switches to very high $\cfl$
numbers, but reduces the timesteps appropriately as soon as a perturbation
enters the flow field.

We combine our time-adaptation with the spatial adaptive multiresolution
technique \cite{Mueller:02}. 
The paper is organized as follows. We start with a brief description of the
fluid equations and their discretization by implicit finite volume schemes, see
Section~\ref{FAMS}.  The adjoint error control is presented in Section
\ref{aec}, including a complete error representation (Section~\ref{sec.error_representations}) and a
related space-time splitting (Section~\ref{sec.space_time_splitting}). Section~\ref{sect:conservative_dual_problem}
contains the conservative approach to the dual problem und its boundary conditions. In Section~\ref{section.numexp} we present the adaptive method in time. 
In Section \ref{section.numexp_setup} we present the instationary test case, a 2D Euler
transonic flow around a circular arc bump in a channel. In
Section~\ref{section.numexp_fully_implicit} results of the fully implicit and a
mixed explicit-implicit time adaptive strategy are presented to illustrate the
efficiency of the scheme.  In Section~\ref{section:conclusion} we summarize
our results.

We refer the reader to \cite{Steiner_2008_diss, SteinerNoelle, SteinerMuellerNoelle}
for further details and references.

{\bf Acknowledgement:} We would like to thank Ralf Hartmann and Mario Ohlberger for stimulating discussions.

\section{Governing equations and finite volume scheme}
\label{FAMS}

We will work in the framework of hyperbolic systems of conservation laws in space and time,
\mm{
  \label{eq:cl}
  U_t + \nabla \cdot f(U) =0 \quad \textnormal{in } \Omega_T.}
Here $\O\subset\R^d$ is the spatial domain with boundary $\Gamma\!:=\!\partial
\O\subset\R^d$ and $\Omega_T = \Omega\times[0,T)\subset\Rd1$ is the space-time
domain with boundary $\Gamma_T\!:=\!\partial\OT\subset\Rd1$. $U$ is the vector of conservative variables
and  $f$ the array of the corresponding convective fluxes $f_i$, $i= 1,\dots,d$,
in the $i$th coordinate direction. Our prime example are the multi-dimensional Euler equations of gas dynamics.

We prescribe boundary conditions on the incoming characteristics as follows:
\mm{
  \label{eq:cl_bc}
  P_-(U^+)\,(\fn(U^+) -g)=0\quad \textnormal{on } \Gamma_T.}
Here $\ff(U):= (f(U),U)$ is the space-time flux, $\nu$ the space-time
outward normal to $\Omega_T$, and $\fn(U) := \ff(U)\cdot\nu$ the space-time
normal flux. $U^+$ is the interior trace of $U$ at the boundary $\GT$ (or any
other interface used later on). 
Given the boundary value $U^+$ and the corresponding Jacobian matrix $\fn'(U^+)$,
let $P_-(U^+)$ be the $(d+2)\times(d+2)$-matrix which realizes the
projection onto the eigenvectors of $\fn'(U^+)$ corresponding to negative
eigenvalues.
Then the matrix-vector product $P_-(U^+)\, \fn(U^+)$ is the incoming component
of the normal flux at the boundary, and it is prescribed in \eqref{eq:cl_bc}.
See \cite{SteinerMuellerNoelle} for details.

We approximate \eqref{eq:cl}--\eqref{eq:cl_bc} by a first or second order finite volume
scheme with implicit Euler time discretization.   The computational spatial grid
$\Oh$ is a set of open cells $V_i$ such that
$$
\bigcup_{i}\;\overline{V_i} = \Ob.
$$
The intersection of the closures of two different cells is either empty or a
union of common faces and vertices. Furthermore let $\myN(i)$ be the set of
cells that have a common face with the cell $i$, $\partial V_i$ the boundary of
the cell $V_i$ and for $j\in\myN(i)$ let $\Gamma_{ij}:= \partial V_i  \cap
\partial V_j $ be the interface between the cells $i$ and $j$ and $n_{ij}$ the
outer spatial normal to $\Gamma_{ij}$ corresponding to cell $i$. Since we will
work on curvilinear grids, we require that the geometric consistency condition 
\mm{\sum_{j \in N ( i )} |\Gamma_{ij}| n_{ij} = 0
\label{eq:geometric_consistency}}
holds for all cells.

Let us define a partition of our time interval $I:=(0,T)$ into subintervals
$\In = [t_\no, t_\nn]$, $1\leq \nn \leq N$, where
\begin{align*}
0=t_0<t_1<\ldots<t_\nn<\ldots<t_N=T.
\end{align*}
The timestep size is denoted by  $\dtn:= t_\nn-t_\no$.
Later on this partition will be defined automatically by the
adaptive algorithm. 
We also denote the space-time cells and faces by $\Vin:=\Vi\times\In$ and 
$\Gijn:=\Gij\times\In$, respectively.
Given this space-time grid the implicit finite volume discretization of
\eqref{eq:cl} can be written as
\mm{
  U_i^\nn+\frac{\dtn}{|V_i|}\,\sum_{j\in\myN(i)}
  |\Gamma_{ij}|\, F_{ij}^\nn
  = U_i^\no \;\;\;\text{for}\;\;\nn\geq 1.
  \label{eq:euler_fv}}
It computes the approximate cell averages $U^\nn_i$ 
of the conserved variables on the new time level. 
For interior faces $\Gij$, the canonical choice for the numerical flux is a
Riemann solver,
\mm{\label{eq:euler_fv_fint} F_{ij}^\nn := F_{riem}(U_{ij}^\nn,U_{ji}^\nn,n_{ij})}
consistent with the normal flux $f_n(U) = f(U)\cdot n_{ij}$.
In the numerical experiments in Section \ref{section.numexp_fully_implicit} we choose Roe's solver \cite{Roe1981}.
If $\Gij\subset\partial\OT=:\GT$, then we follow the definition of a weak solution
and define the numerical flux at the boundary by
\begin{align}\label{eq:euler_fv_fext}
  F_{ij}^\nn := P_+(U_{ij}^\nn)\,\ff_{\nu_{ij}}(U_{ij}^\nn)
  + P_-(U_{ij}^\nn)\,g_{ij}^\nn,
\end{align}
where $g_{ij}^\nn$ is the average of $g$ over $\Gijn$. 

For simplicity of presentation we neglect in our notation that due to higher
order reconstruction the numerical flux usually depends on an enlarged stencil
of cell averages.

\section{Adjoint error control - adaptation in time} \label{aec}

In order to adapt the timestep sizes we use a method which involves adjoint
error techniques. We have applied this approach successfully to Burgers'
equation in \cite{SteinerNoelle}, and to the Euler equations of gas dynamics in \cite{SteinerMuellerNoelle}.

Since a finite volume discretization in space and a backward Euler step in time
are a special case of a Discontinuous Galerkin discretization, techniques based
on a variational formulation can be transferred to finite volume methods.

The key tool for the time adaptive method is a space-time splitting of adjoint
error representations for target functionals due to S\"uli \cite{Sueli} and
Hartmann \cite{Hartmann1998}. It provides an efficient choice of timesteps for
implicit computations of weakly instationary flows. The timestep will be very
large in time regions of stationary flow, and become small when a perturbation
enters the flow field.

\subsection{Variational Formulation}

In this section we rewrite the finite volume method as a Galerkin method, 
which makes it easier to apply the adjoint error control techniques.

Let us first introduce the space-time numerical fluxes. Let
$\Oin=V_i\times\In\in\OTh$ be a space-time cell, and let
$\gamma\subset\partial\Oin$ be one of its faces, with outward unit normal $\nu$.
There are two cases: if $\nu$ points into the spatial direction, then
$\gamma=\Gij\times\In$ and $\nu=(n,0)$. If it points into the positive time
direction, then $\gamma=V_i\times\{t_\nn\}$, and $\nu=(0,1)$. Now we define the
space-time flux by
\mm{
  \FFn(U_h) = \left\{ \begin{array}{cl}
  F_{ij}^\nn \;\; \text{from} \;\; \eqref{eq:euler_fv_fint} & 
  ~~\text{if}~~\nu=\nij \;\; \text{and}~~\gamma\in\Eint \\
  F_{ij}^\nn \;\; \text{from} \;\; \eqref{eq:euler_fv_fext} & 
  ~~\text{if}~~\nu=\nij \;\; \text{and}~~\gamma\in\Eext \\
  (1-\theta)U_i^\no + \theta U_i^\nn & 
  ~~\text{if}~~\nu=(0,1) \;\; \text{and}~~ \nn \geq 1 \\
  U_i^0 & ~~\text{if}~~\nu=(0,1) \;\; \text{and}~~ \nn=0
  \end{array}\right.
  \label{eq:space-time-flux}}
where $\Eint$ are the interior faces and $\Eext$ the boundary faces. In the
third case, $\theta\in[0,1]$, so the numerical flux in time direction is a
convex combination of the cell averages at the beginning and the end of the
timestep. Different values of $\theta$ will yield different time
discretizations, e.g. explicit Euler for $\theta =0$, implicit Euler for $\theta
=1$.

Let $\myVVh:=W^{1,\infty}(\OTh)$ be the space of piecewise Lipschitz-continuous
functions. Now we introduce the semi-linear form $N$ by
\mm{
  N &: \;\;\;\;\;\; \UU \times \myVVh \to \mathbf R \nonumber \\
  N(U,\vp) &:= \sum_\inn(\FFn(U),\vp)_{\Gin} \nonumber \\
  & \;\; - \sum_\inn\left( (U,\vp_{h,t})_{\Oin} + 
  (f(U),\nabla\vp)_{\Oin}. \right)
  \label{eq:var_form}}

Here and below the sum is over the set $\{(i,\nn)\,|\,\Oin\in\OTh\}$, i.e. all
gridcells. Now we rewrite the finite volume method \eqref{eq:euler_fv} as a
first order Discontinuous Galerkin method (DG0): 

\mm{
  \text{Find}\;\; U_h\in\VVh0 \;\;\;\;\text{such that}\;\;\;\;
  N(U_h,\vp_h) = 0  \quad\quad\forall \vp_h\in\VVh0,
  \label{eq:DG0}}
where $\VVh0$ is the space of piecewise constant functions over $\OTh$.
\begin{remark}
For the DG0 method, $U_h,\vp\in\VVh0$ are piecewise constant, so the last two
terms in \eqref{eq:DG0}, containing derivatives of $\vp$, disappear. Moreover,
due to the geometric condition \eqref{eq:geometric_consistency}
\mmn{
\sum\limits_{\{j\,|\,\Gijn\subset\Gin\}}\ff(U_i^\nn)\cdot\nu_{ij}^\nn = 0}
holds for all cells $\Oin$. Therefore, the DG0 solution may be characterized by:
Find $U_h\in\VVh0$ such that
\mm{
  \sum_\inn(\FFn(U)-\fn(U^+),\vp_h)_{\Gin} = 0
  \quad\quad\forall \vp_h\in\VVh0. \label{eq:DG0_b}}
This form is convenient to localize our error representation later on.  
\end{remark}

\subsection{Adjoint error representation for target functionals}
\label{sec.error_representations}

In this section we define the class of target functionals $J(U)$ treated in
this paper, state the corresponding adjoint problem and prove the
error representation which we will use later for adaptive timestep control.

Before we derive the main theorems, we would like to give a preview of an
important difference between error representations for linear and nonlinear
hyperbolic conservation laws.
For linear conservation laws (and many other linear PDE's), it is possible to
express the error in a user specified functional,
\mm{\label{eq:error_representation_a}
  \eJ &:= J(U)-J(U_h),}
as a computable quantity $\eta$, so
\mm{\label{eq:eJ}
  \eJ &= \eta}
(see e.g. \cite{BeckerRannacher1996, SteinerMuellerNoelle, SueliHouston2003} and the references therein).
In general, $\eta$ will be an inner product of the numerical residual
with the solution of an adjoint problem. Below we will see that such a 
representation does not hold for nonlinear hyperbolic conservation laws.
The nonlinearity will give rise to an additional error $\egm$ on the
inflow boundary, an error $\egp$ on the outflow boundary, and a linearization
error $\eO$ in the interior domain:
\begin{theorem}\label{theorem:error_representation_h}
Suppose $\vp \in \VV$ solves the approximate adjoint problem
\mm{
  \partial_t \vp + \At^T\nabla\vp &= \psi \quad \; \; \mathrm{in} \; \OT,
  \label{eq:adjoint_h} \\
  \Pt^T_+ (\vp -\psi_\Gamma) &=0  \quad \; \;\mathrm{on} \; \Gamma_T,
  \label{eq:adjoint_h_bc}}
where $\At := A(U_h)$. Let
\mmn{
  \egm &:= -((P_-(U^+)-P_-(U_h^+))g,\vp)_\GT \\
  \egp &:= -(P_+(U^+)\fn(U^+)-P_+(U_h^+)\fn(U_h^+),\vp-\psi_\Gamma)_\GT \\
  \eO &:= (f(U)-f(U_h)-\At(U-U_h),\nabla\vp)_\OT}
and
\mm{\eta &:= N(U_h,\vp) = N(U_h,\vp-\vp_h). \label{eq:eta}} 
Then
\mm{\label{eq:error_representation_h}
  \eJ+\egm+\egp+\eO = \eta.}
\end{theorem}
Our adaptation is based on computing and equidistributing this $\eta$.
Typical examples for the functional $J$ are the lift or the drag of a
body immersed into a fluid. To simplify matters we consider
functionals of the following form:
\mm{
  J(U) = (U,\psi)_{\OT} - (P_+(U^+) \fn(U^+),\psi_\Gamma)_{\GT},
  \label{eq:functional}}
where $\psi$ and $\psi_\Gamma$ are weighting functions in the interior
of the space-time domain $\Omega_T$ and at the boundary $\Gamma_T$.
For the proof of Theorem~\ref{theorem:error_representation_h} as well as an illustrative example of the functional $J$ we refer again to \cite{SteinerMuellerNoelle}.

For the adjoint problem \eqref{eq:adjoint_h} and \eqref{eq:adjoint_h_bc} the role
of time is reversed and hence $\Pt_+$ plays the role of $P_-$ in
\eqref{eq:cl_bc}. Here $\psi_\Gamma$ comes from the weighting function in
the functional \eqref{eq:functional}.

\subsection{Space-time splitting}
\label{sec.space_time_splitting}

The error representation \eqref{eq:error_representation_h} is not yet suitable for
time adaptivity, since it combines space and time components of the residual and
of the difference $\vp-\vp_h$ of the dual solution and the test function. The
main result of this section is an error estimate whose components depend either
on the spatial grid size $h$ or the timestep $k$, but never on both. The key
ingredient is a space-time splitting of \eqref{eq:error_representation_h} based on
$L^2$-projections. Similar space-time projections were introduced previously in
\cite{Hartmann1998,Sueli}. In \cite{SteinerNoelle} we adapted them to the finite
element spaces and space-time Discontinuous Galerkin methods of arbitrary order.

The splitting takes the form
\mm{\label{eq:error_representation_hk}
  \eta = \eta_k+\eta_h.}
For brevity, we only present the details for first-order finite volume schemes for which we obtain  
\mm{
  \eta_k &= \frac{\dtn}2\,\sum_\inn
  \left((1-\theta)(U_i^\no-U_i^{\nn-2})
  +\theta(U_i^\nn-U_i^\no),\psi-\At^Tw\right)_{V_i}
  \label{eq:eta_k}}
(we set $U_i^{-1}=U_i^{0}$ in the first summand). Thus our temporal error
indicator is simply a weighted sum of time-differences of the approximate
solution $U_h$, and the weights can be computed from the data $\psi$ and the
solution $w$ of the conservative dual problem~\eqref{eq:adjoint_w}. 

In our adaptive strategy, we will use the localized indicators
\mm{
  \ebkn := \frac12\,\sum\limits_i
  \left|((1-\theta)(U_i^\no-U_i^{\nn-2})
  +\theta(U_i^\nn-U_i^\no),\psi-\At^Tw)_{V_i}\right|.
  \label{eq:localized_error_indicator}}
In the next section we present an example for the boundary conditions for the
dual problem.

\section{The conservative dual problem}
\label{sect:conservative_dual_problem}

In this section we present the conservative approach to the dual problem, which we introduced in
\cite{SteinerNoelle} and derive boundary conditions for the  gradient of the
dual problem.

The adjoint equation \eqref{eq:adjoint_h} is a system of linear transport
equations with discontinuous coefficients. Therefore, numerical approximations
may easily become unstable. Another inconvenience is that in order to obtain a
meaningful error representation in \eqref{eq:error_representation_h}, the
approximate adjoint solution $\vp$ should not be contained in $\VVh0$.
Therefore, $\vp$ is often computed in the more costly space $\VV_h^1$.

In \cite{SteinerNoelle} we have proposed a simple alternative which helps to
avoid both difficulties. Instead of computing the dual solution $\vp$ we will
compute its gradient
\mmn{w := \nabla\vp,}
which is the solution of the conservative dual problem
\mm{
  w_t+\nabla(\At^Tw) = \nabla \psi \;\;\;\text{in}\;\OT.
  \label{eq:adjoint_w}}
This system is in conservation form, and therefore it can be solved by any
finite volume or Discontinuous Galerkin scheme. Moreover, \eqref{eq:adjoint_w}
may be solved in $\VVh0$, since a piecewise constant solution $w$ already
contains crucial information on the gradient of $\vp$. 

The scalar problem treated in \cite{SteinerNoelle} was set up in such a way that
the characteristic boundary conditions for the dual problem became trivial. In
\cite{SteinerMuellerNoelle} we developed boundary conditions for the more general
initial boundary value problem \eqref{eq:cl} -- \eqref{eq:cl_bc}.

Denoting the flux in \eqref{eq:adjoint_w} by $H:=\At^Tw$, the boundary condition
\eqref{eq:adjoint_h_bc} becomes
\mm{
  \Pt_+^T\,(H-H_\Gamma) = 0 \;\;\;\text{on}\;\GT,
  \label{eq:adjoint_w_bc}}
i.e. we prescribe the incoming component $\Pt_+^TH$. Here $H_\Gamma$ is a given
real-valued vector function, which depends on $\psi_\Gamma$.  However, this
characteristic boundary condition needs to be interpreted carefully. Using
\eqref{eq:adjoint_w} and denoting the interior trace at the
flux by $H_{int}$, we may introduce the boundary flux by \mmn{H := \Pt_-^T \,
H_{int} + \Pt_+^T \, H_{\Gamma} \quad \textnormal{on } \GT.} Note that all the
projections $\Pt_\pm$ used below depend on the point $(x,t)\in\GT$ via
the outside normal vector $\nu(x,t)$. The value $\Pt_-^T H_{int}$ may be assigned
from the trace $w_{int}$ at the interior of the computational domain, $$ \Pt_-^T
\, H_{int} =\Pt_-^T \,(\At^Tw)_{int}. $$ The boundary values $\Pt_+^T \, H_\Gamma$
are computed using the PDE 
\mm{
  \vp_t = -H +\psi
  \label{eq:v_t}}
with boundary values \eqref{eq:adjoint_h_bc},
\mm{
  \Pt_+^T \, H_\Gamma 
  &= \Pt_+^T \,(-\vp_t +
    \psi)|_\Gamma \nonumber \\
  & = -(\Pt_+^T\,\psi_\Gamma)_t 
    + \Pt_+^T\,\psi \nonumber \\
  & \approx - \frac1{\dtn} \left(\Pt_+^{T,\nn}
    \,\psi_\Gamma^\nn - \Pt_+^{T,\no}
    \,\psi_\Gamma^\no\right) 
    + \Pt_+^{T,\no}\,\psi.}
This completes the definition of
the numerical boundary conditions for the conservative dual problem.

\section{Adaptive concept}
\label{section.numexp}

Now we combine the multiscale approach in space \cite{Mueller:02}
and the time adaptive method derived from the space-time splitting of the error
representation to get a space-time adaptive algorithm:
\begin{itemize}
\item solve the primal problem \eqref{eq:cl} on a {\em coarse}
  adaptive spatial grid (e.g. level $L=2$) using uniform $\cfl$ numbers ($\cfl = 0.8$),
\item compute the dual problem \eqref{eq:adjoint_h} and
  \eqref{eq:adjoint_h_bc} and the space-time-error representation 
  \eqref{eq:error_representation_hk}. In particular, compute the localized 
  error indicators $\ebkn$ using \eqref{eq:localized_error_indicator}.
\item compute the new adaptive timestep sizes depending on the temporal part
  of the error representation and the $\cfl$ number on the new grid, 
  aiming at an equidistribution of the error,
\item solve the primal problem using the new timestep sizes on a 
  {\em finer} spatial grid (e.g. level $L=2$).
\end{itemize}

The advantage is, that the first computations of the primal problems and the
dual problem are done on a coarse spatial grid, and therefore have low cost.
These computations provide an initial guess of the timesteps for the computation
on the finer spatial grid. We will restrict the timestep size from below to
$\cfl=0.8$, since smaller timestep sizes only add numerical diffusion to the
scheme and increase the computational cost. Note that all physical effects
already have to be roughly resolved on the coarse grid in order to determine a
reliable guess for the timesteps on the fine grid.

We will deal with some aspects in detail in the numerical examples in Section
\ref{section.numexp_fully_implicit}.

\subsection{Asymptotic decay rates}
\label{sec.SEP}
Since the adaptive strategy outlined in Section \ref{section.numexp}
above depends on assumptions on the assymptotic behavior of the error, we first 
try to estimate these decay rates. There is no analytical result which shows how the error terms
$\bar\eta_k$ and $\bar\eta_h$ depend on $k$ and $h$. Therefore, we estimate this 
dependence numerically. We compute a perturbed shock of Burgers equation, for details see \cite{SteinerNoelle}.
We compare the two approaches:
\begin{itemize}
\item refinement only time 
\item and refinement only space.
\end{itemize}

\begin{figure}
  \begin{center}\vspace{-3cm}
  \includegraphics[width=0.49\textwidth]{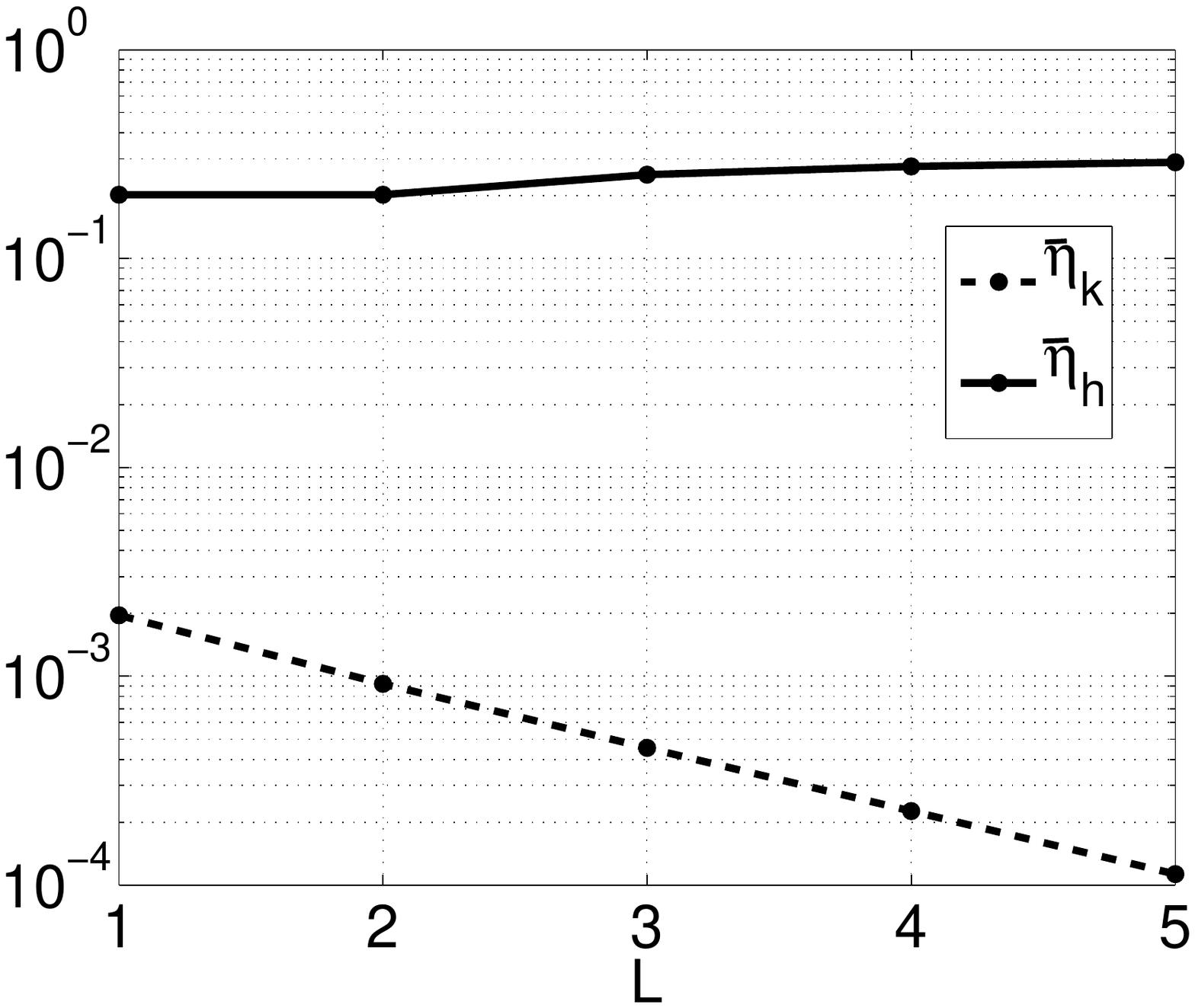}
  \includegraphics[width=0.49\textwidth]{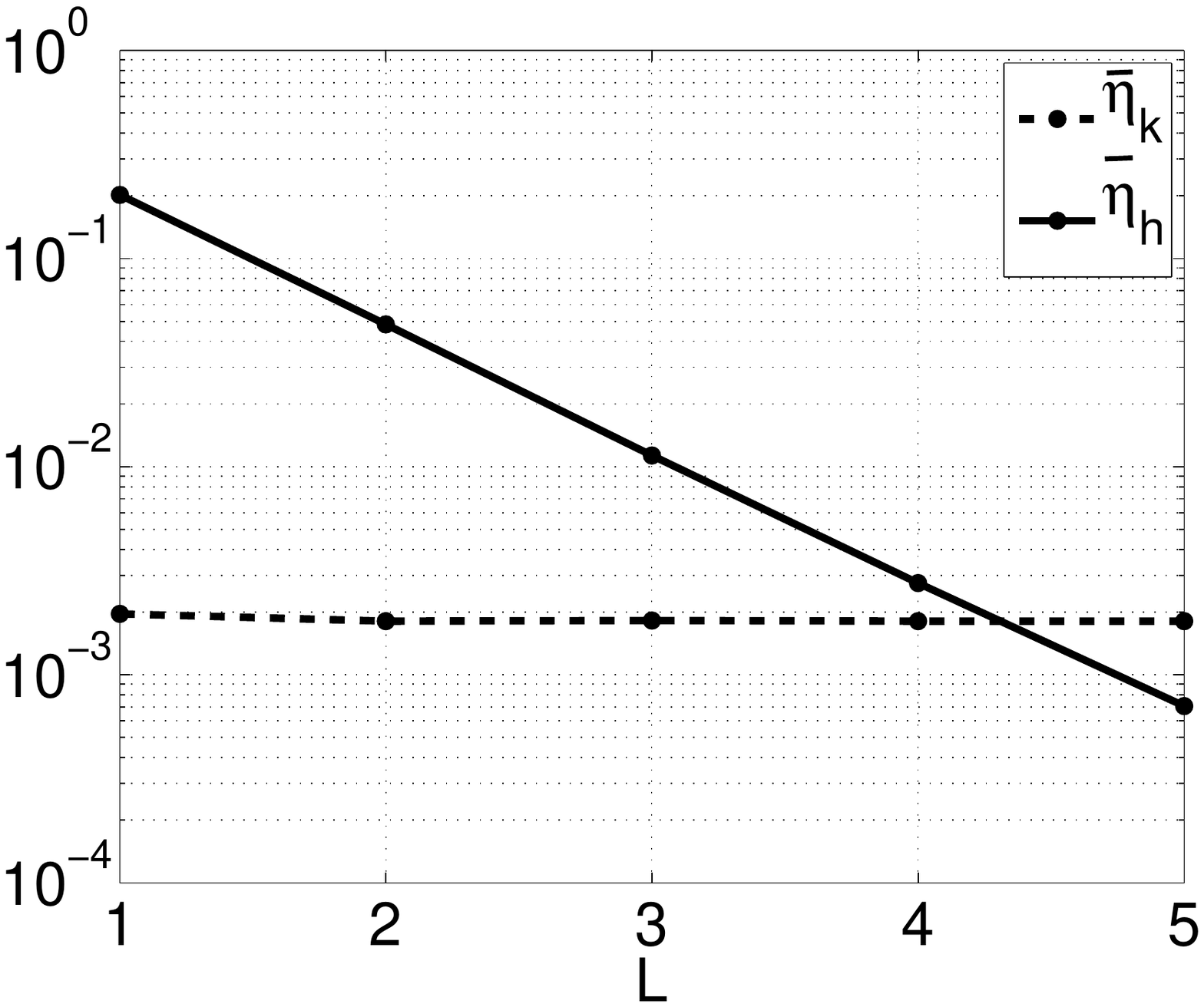}\vspace{-2cm}
  \caption{Error representation for Burgers equation, 
first order method, $\bar{\eta_k}$ and $\bar{\eta_h}$ versus level 
 of refinement. Left: uniform refinement in time. Right: uniform refinement in space \label{fig.burgers_1_etabar}}
  \end{center}
\end{figure}

\begin{table}  
\begin{center}
\begin{tabular}{c|c|c|c|c|c|c|c}
$L$ & $\bar{\eta_k}$	&$\bar{\eta_h}$&$\eta_k$	&$\eta_h$	&$
J(u_h)$ &$\eta_h+\eta_k$& $\theta$ \\ \hline
1&     1.96e-03&  2.02e-01&  1.29e-04&  2.00e-01&  1.72e+00&  2.01e-01&  5.57e+00   \\  \hline 
2&     9.81e-04&  4.83e-02&  1.83e-05&  4.75e-02&  1.74e+00&  4.75e-02&  5.49e+00   \\  \hline 
3&     4.81e-04&  1.21e-02&  3.57e-06&  1.17e-02&  1.75e+00&  1.17e-02&  5.71e+00   \\  \hline 
4&     2.37e-04&  3.10e-03&  1.30e-06&  2.89e-03&  1.75e+00&  2.89e-03&  7.06e+00
\end{tabular}\\ 
\caption{Efficiency $\theta = \frac{\eta_h+\eta_k}{J(u)-J(u_h)}$of the error
representation, $\cfl = 0.8$. } 
\end{center}
\end{table}

Each of the plots in Figure~\ref{fig.burgers_1_etabar} show the error
estimators $\bar\eta_k$ (error in time) and  $\bar\eta_h$ (error in space).  In the
Figure~\ref{fig.burgers_1_etabar}(a) we refined only in time. Here the spatial
error remains constant, while the time error still decreases with first order.
The second Figure \ref{fig.burgers_1_etabar}(b) shows the refinement only in space. 
The time error $\bar\eta_k$ is almost constant, while the spatial error is decreasing with second order. 

Numerically the terms $\bar\eta_t$  and $\bar\eta_h$ behave as expected. 
They depend either on $k$ or on $h$, but never on both. The behaviour of
$\eta_h$ and $\eta_k$ is very similar, and not displayed here.

\begin{remark}
The numerically validated results can be used for adaptive 
grid refinement. The error estimator $\bar\eta_h$ can be used as an indicator for spatial
adaption and the estimator $\bar\eta_k$ for time step control.
\end{remark}

\section{Setup of the numerical experiment}
\label{section.numexp_setup}

An instationary variant of a classical stationary 2D Euler transonic flow \cite{RV:81}, is investigated to illustrate the efficiency of the
adaptive method. 

\begin{figure}[ht]
    \begin{center}
          \includegraphics[width=10cm,clip]{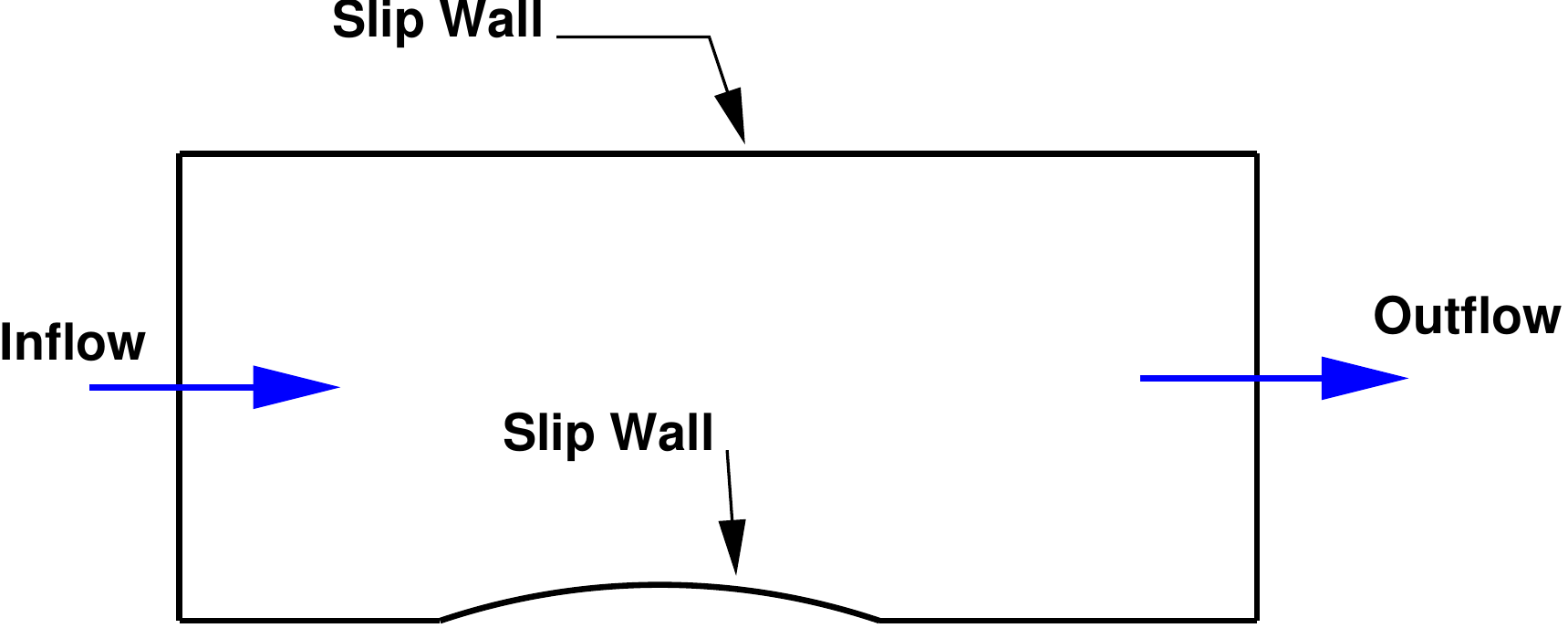}
    \end{center}
  \caption{Circular arc bump configuration of the computational domain $\Omega$. \label{channel-config}}
\end{figure}

{\bf Steady state configuration.} First we consider the classical setup in the
stationary case. The computational domain is a channel of $3m$ length and $2m$
height with an arc bump of $l = 1m$ secant length and $h = 0.024m$ height cut
out, see Figure \ref{channel-config}.  At the inflow boundary, the Mach number
is 0.85 and a homogeneous flow field characterized by the free-stream quantities
is imposed. At the outflow boundary, characteristic boundary conditions are
used. We apply slip boundary conditions across the solid walls, i.e., the normal
velocity is set to zero. In the numerical examples in Section
\ref{section.numexp_fully_implicit} the height of the channel is $2m$ and the
length $6m$.

\begin{figure}[ht]
    \begin{center}\vspace{-5cm}
          \includegraphics[width=\textwidth]{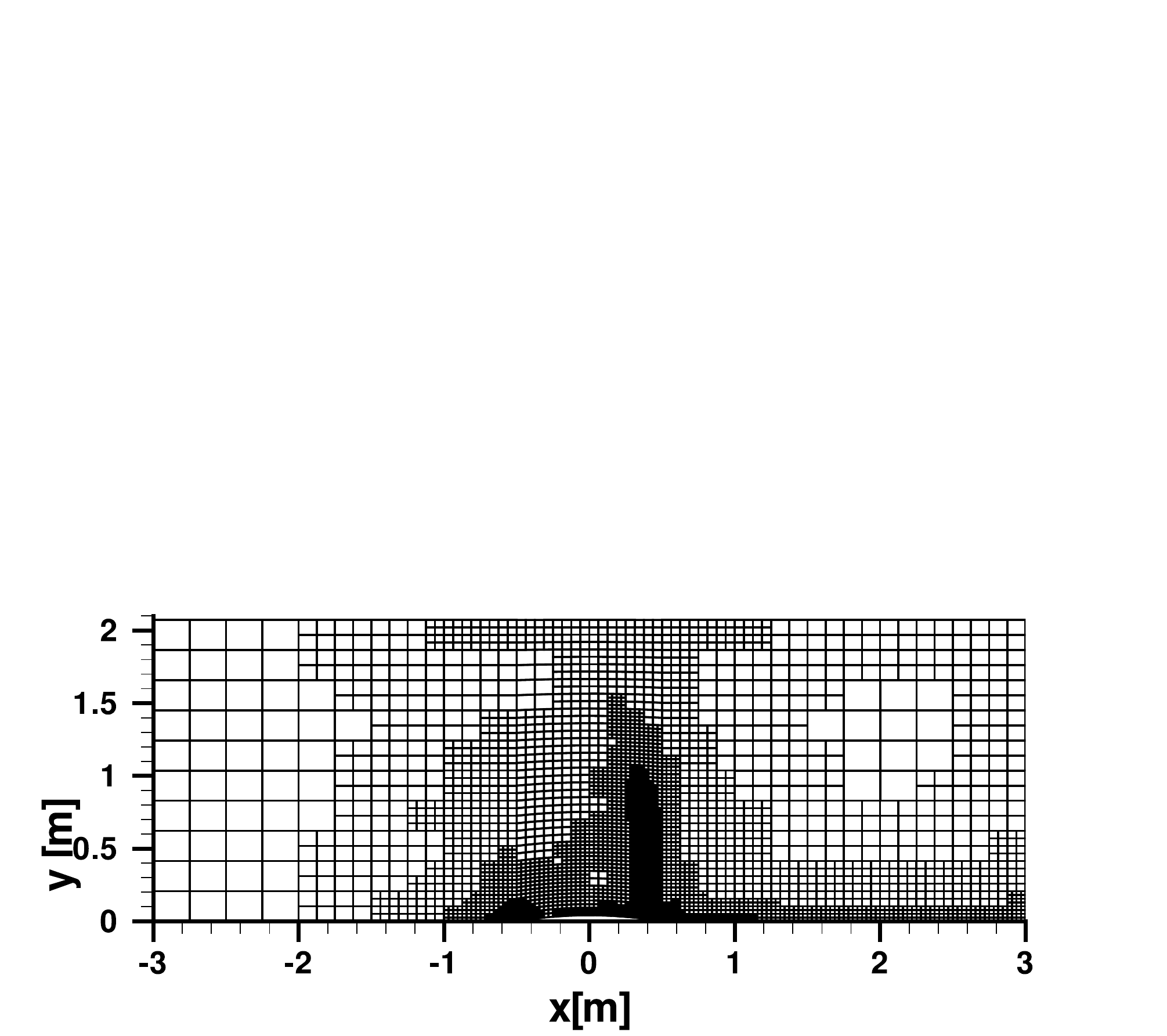}
    \end{center}
  \caption{Adaptive grid $L=5$, to steady state solution in Figure \ref{channel-stat1} of the circular arc bump configuration.  \label{channel-statgrid1}}
\end{figure}

\begin{figure}[ht]
    \begin{center}\vspace{-5cm}
          \includegraphics[width=\textwidth]{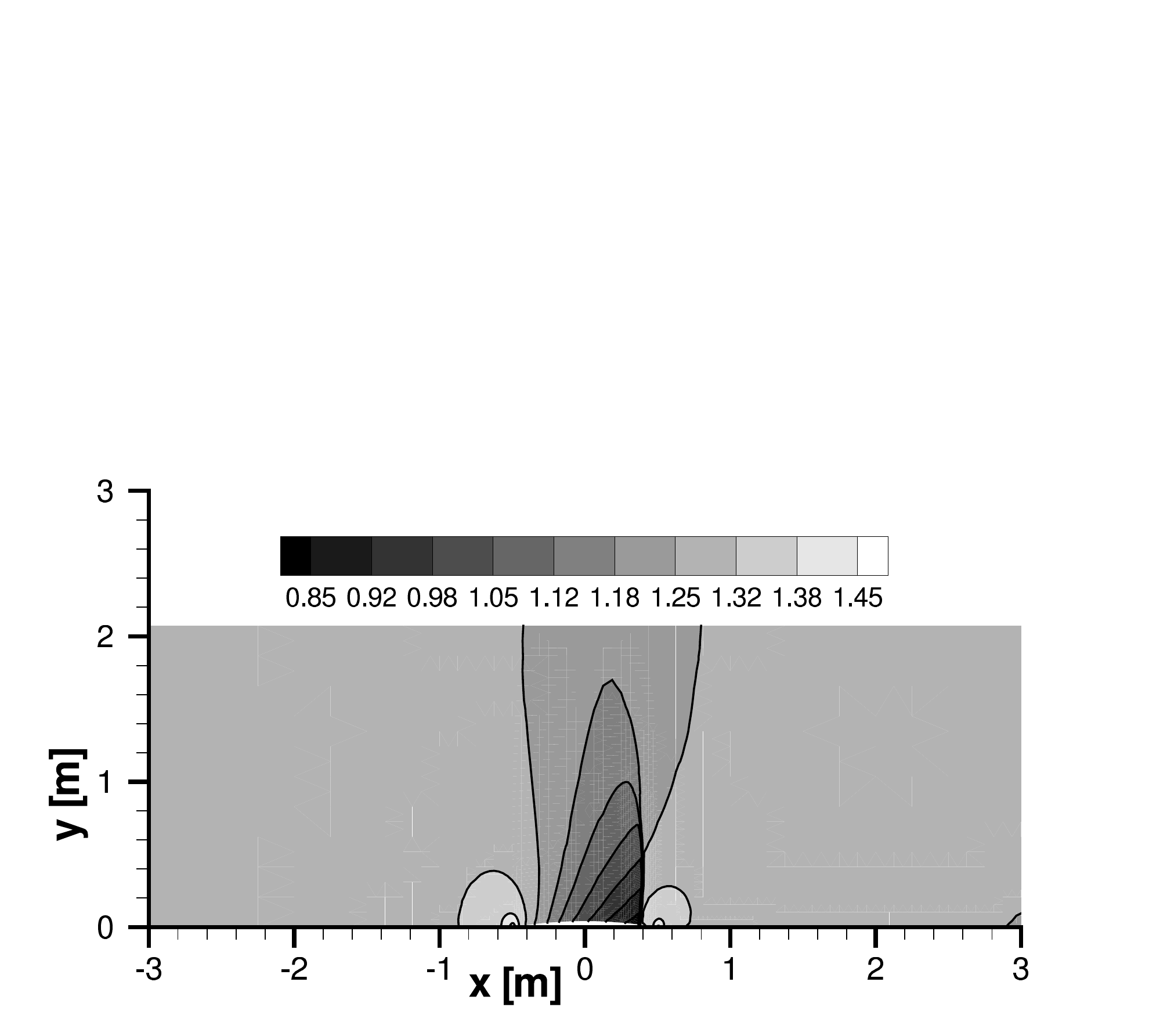}
    \end{center}
  \caption{Steady state solution of the circular arc bump configuration: Isolines of the density, $L=5$.  \label{channel-stat1}}
\end{figure}

The threshold value in the grid adaptation step for the multiscale analysis  is
$\varepsilon = 1 \times 10^{-3}$ and computations are done on adaptive grids
with finest level $L = 2$ and $L=5$ respectively. In general, a smaller
threshold value results in more grid refinement whereas a larger value gives
locally coarser grids. 

In the stationary case at Mach 0.85 there is a compression shock separating a
supersonic and  a subsonic domain. The shock wave is sharply captured  and
the stagnation areas are highly resolved, see Figures \ref{channel-statgrid1} and
\ref{channel-stat1}.  

We will use this steady state solution as initial data for the instationary test
case.

\begin{figure}
  \begin{center}
  \includegraphics[width=\textwidth]{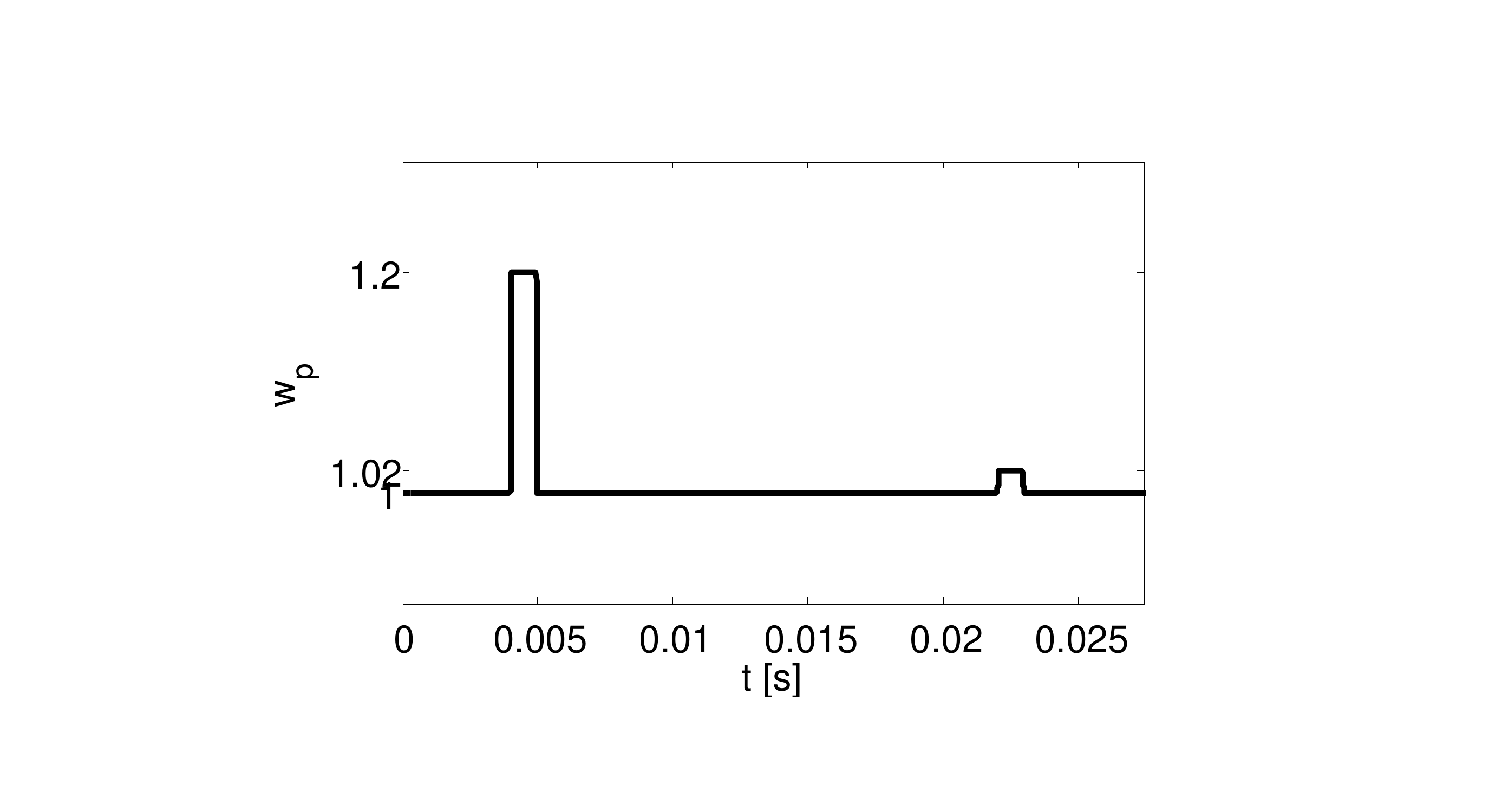}
  \vspace*{-15mm}
  \caption{Weighting function of the perturbation. $p_{in}(t) = p_\infty w_p(t)$} \label{fig.pertub_p}
  \end{center}
\end{figure}

\begin{figure}[ht]
  \begin{center}\vspace*{-30mm}
  \includegraphics[width=0.6\textwidth]{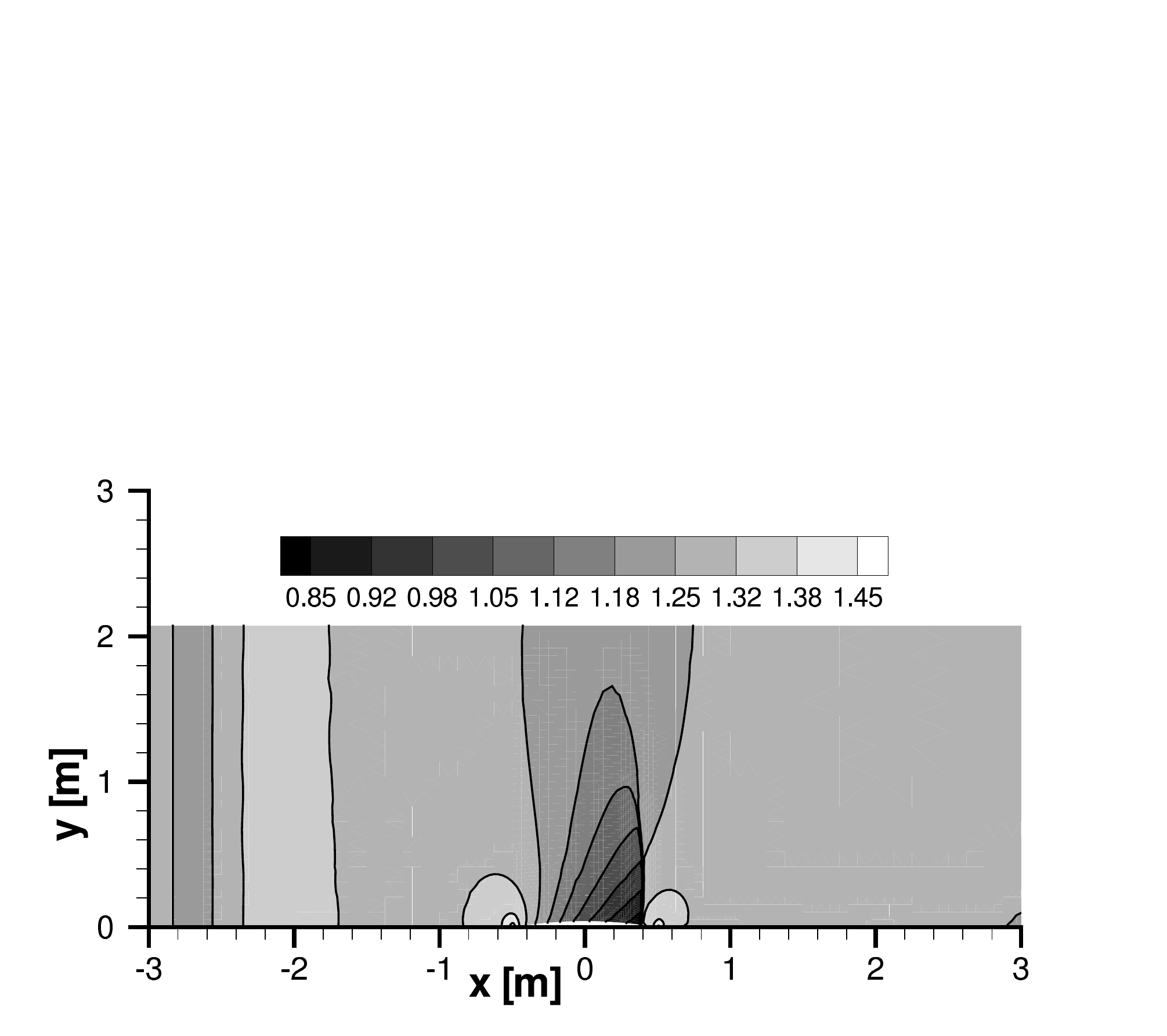}\vspace*{-30mm}
  \includegraphics[width=0.58\textwidth]{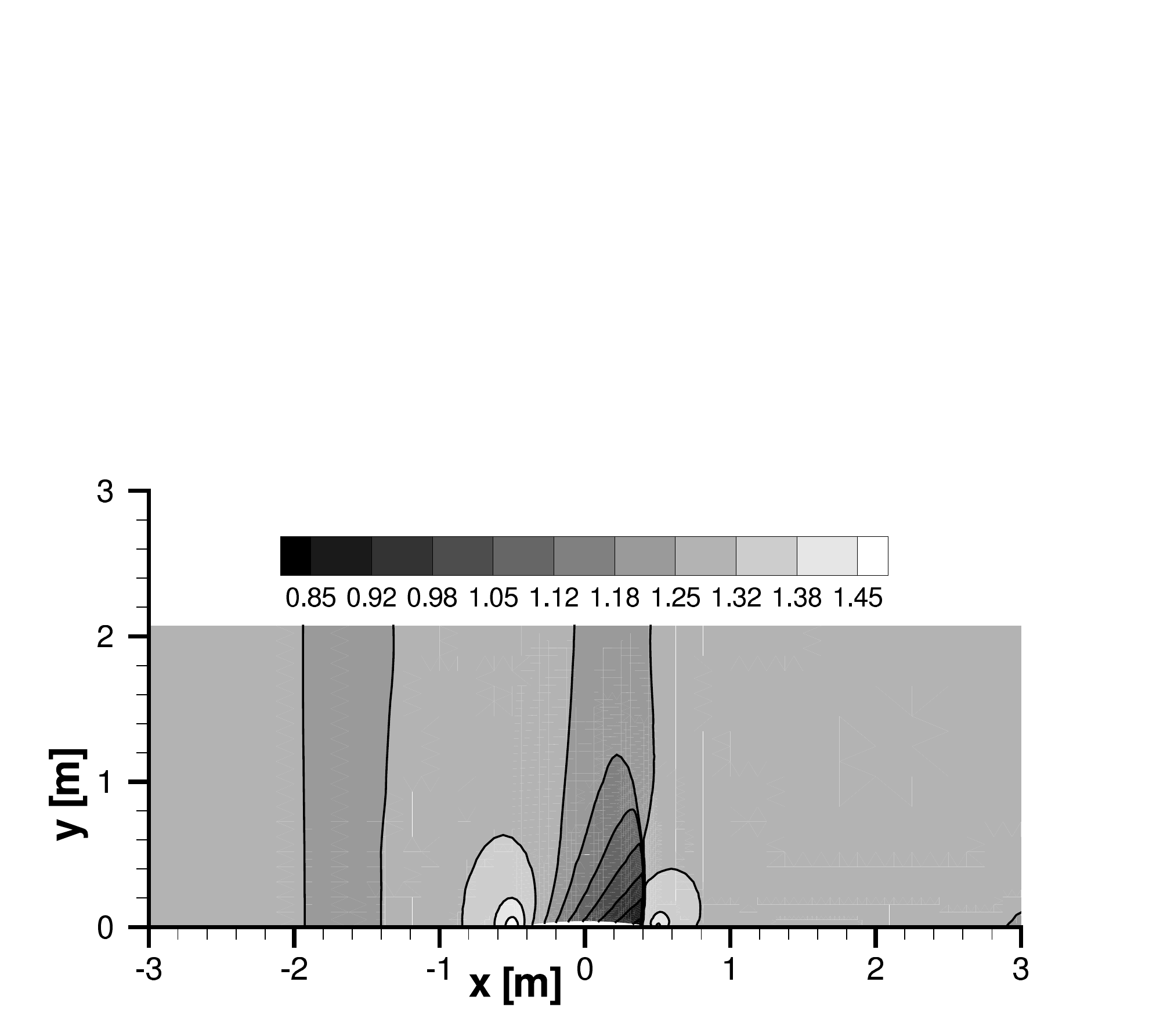}\vspace*{-30mm}
  \includegraphics[width=0.58\textwidth]{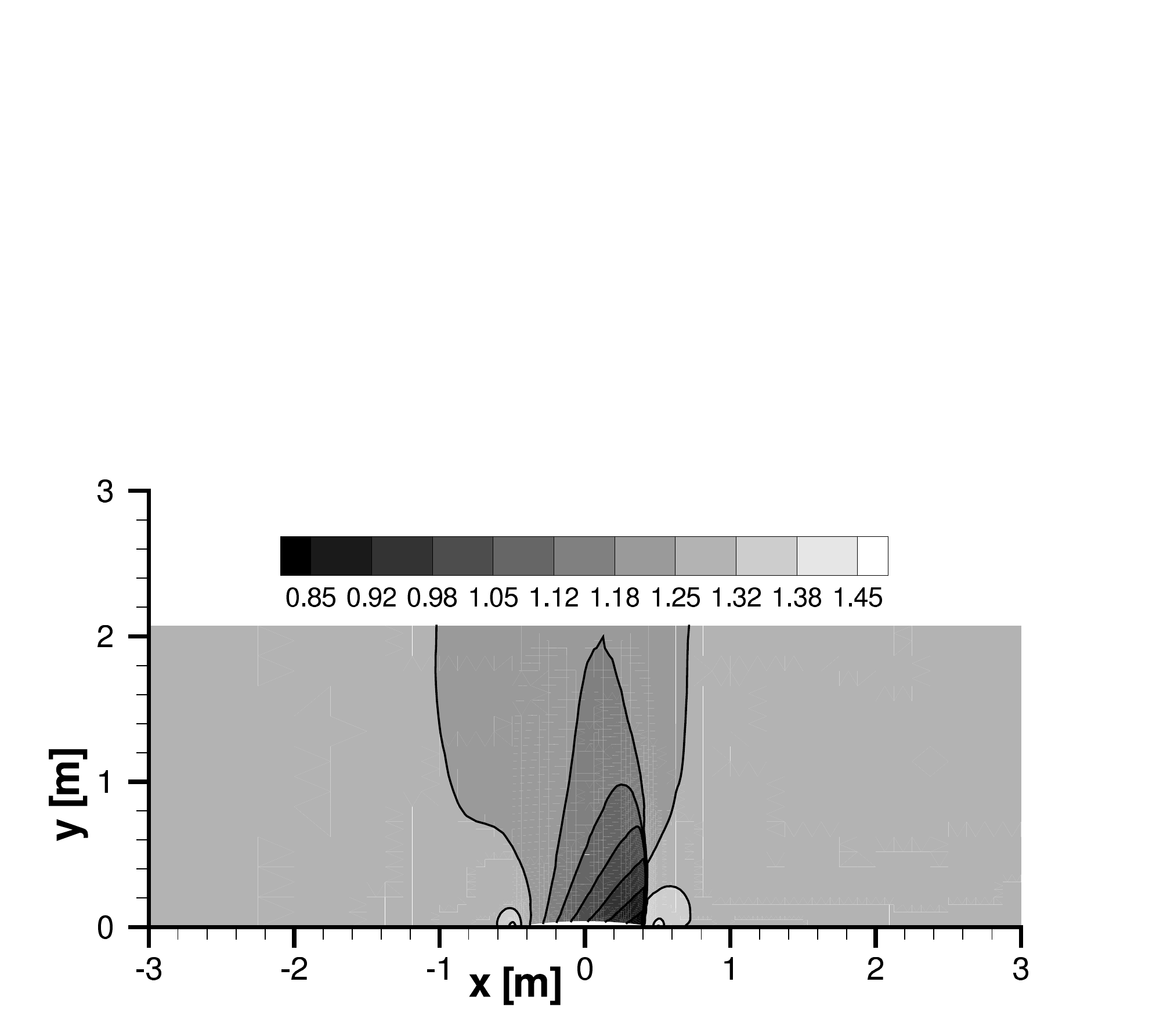}\vspace*{-30mm}
  \includegraphics[width=0.58\textwidth]{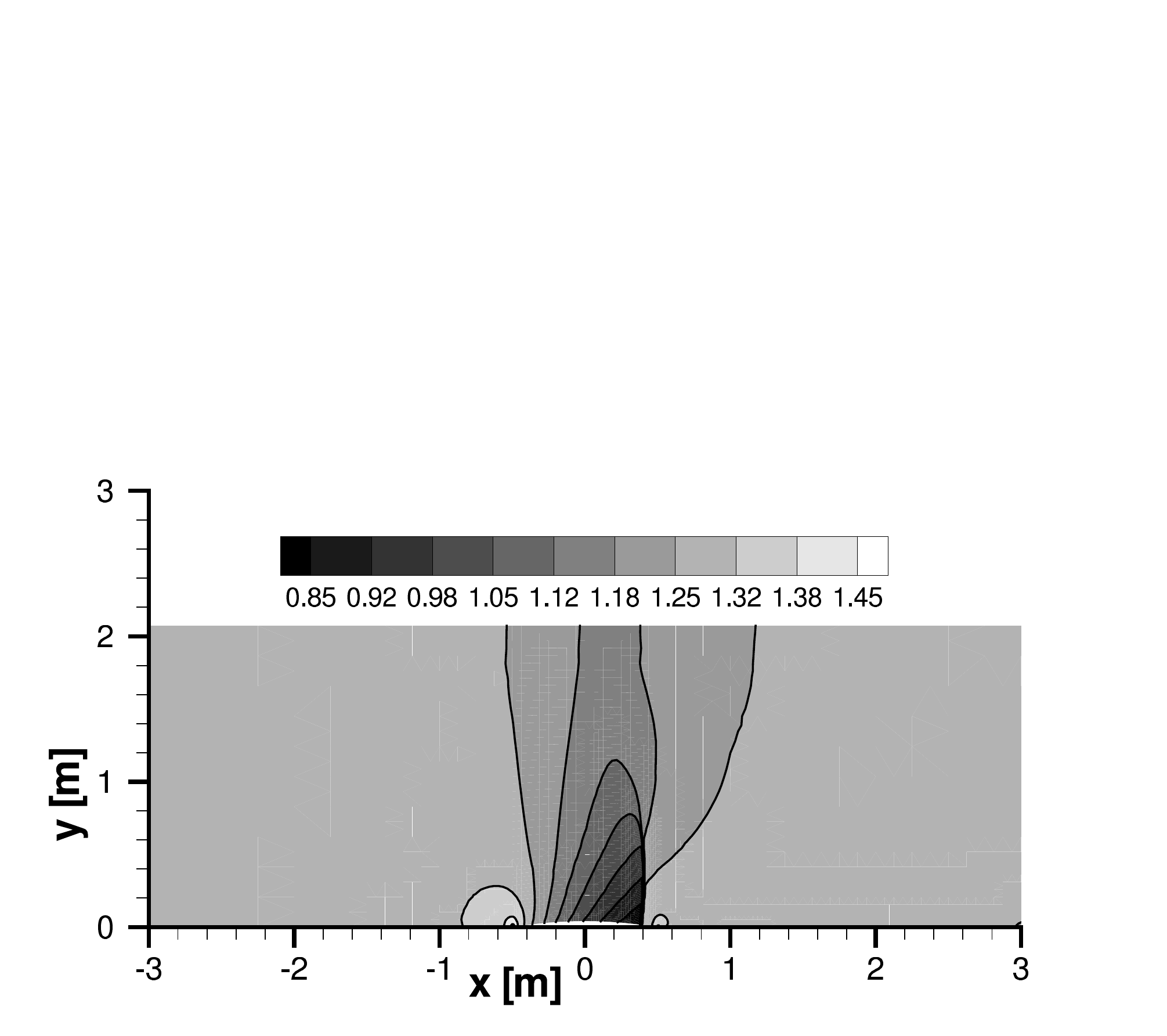}\vspace*{-30mm}
  \includegraphics[width=0.58\textwidth]{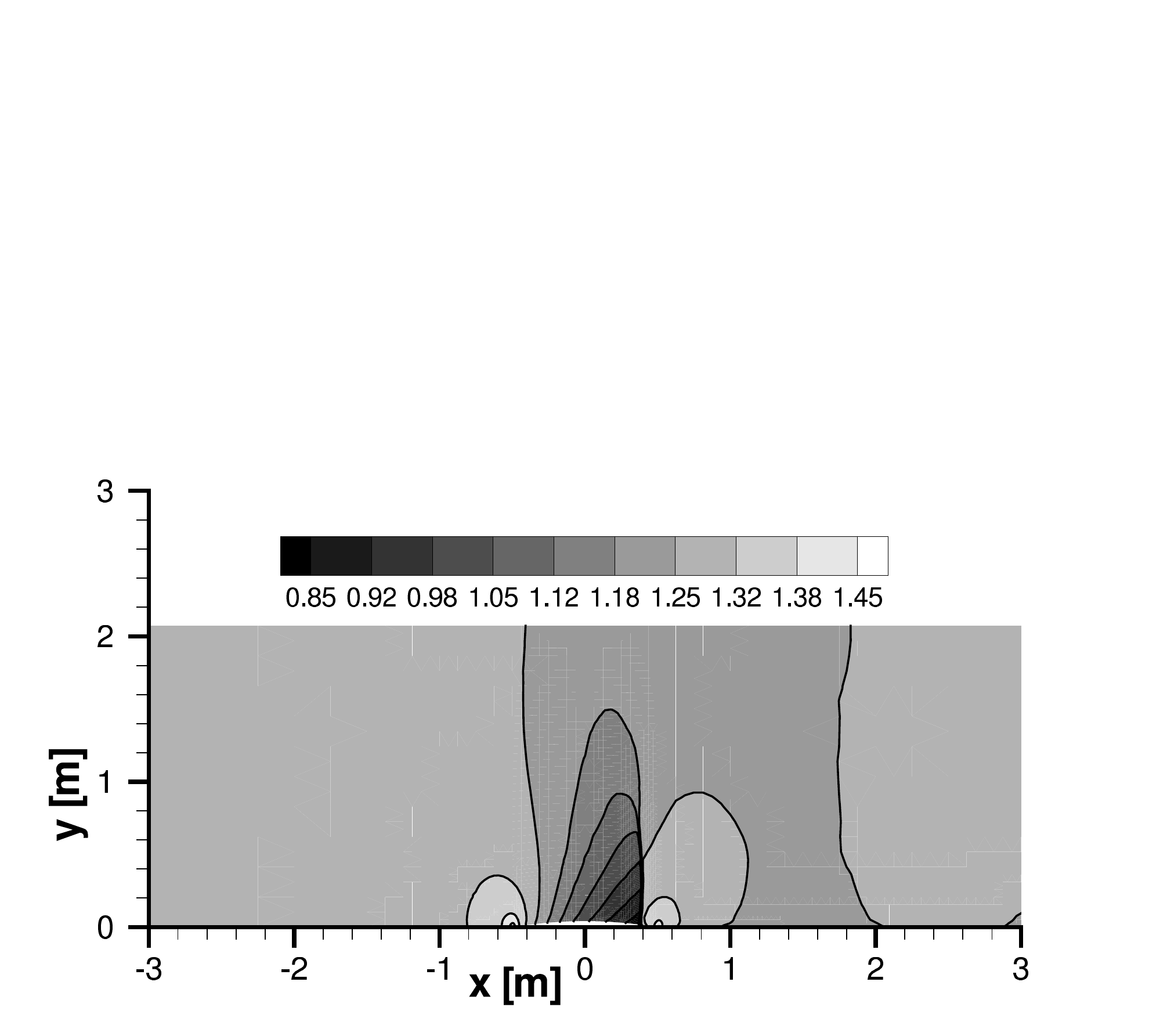}\vspace*{-5mm}
  \end{center}
  \caption{Instationary solution of the circular arc bump configuration, uniform timestep $\cfl =1$, adaptive spatial grid on level $L=5$, isolines of the density, perturbation entering on the left and leaving on the right side of the computational domain, from top to bottom: $t$ = 0.0057s, 0.00912s, 0.01254s, 0.01596s, 0.01938s.}\label{channel-stat3}
\end{figure}

{\bf Instationary test case.} Now we define our instationary test case
prescribing a time-dependent perturbation coming in at the inflow boundary.
First we keep the boundary conditions fixed, and prescribe the corresponding
stationary solution as initial data. Then we introduce for short time periods perturbations of the pressure at the left boundary. 
The first perturbation is about 20 percent of
the pressure at the inflow boundary and the second 2 percent, see Figure \ref{fig.pertub_p}. The perturbations
imposed move through the domain and leave it at the right boundary. Then the
solution is stationary again, see Figure \ref{channel-stat3}. The total time is $t= 0.029s$. The formulas of the pertubations are given in detail in \cite{SteinerMuellerNoelle}.

The first computation is done on an adaptive grid with finest level $L=2$. We
also compute the dual solution and the error representation on this level. Using
the time-space-split error representation \eqref{eq:error_representation_hk} we
derive a new timestep distribution aiming at an equidistribution of the error.
Finally this is modified by imposing a $\cfl$ restriction from below.

We aim to equidistribute the error and prescribe a tolerance
$\;\tol(5)=2^{-3}\bar\eta_k^{ref}\,$, where $\bar\eta_k^{ref}$ is the temporal
error from the computation on level $L=2$.

For this set-up we will show that the adaptive spatial refinement together with
the time-adaptive method will lead to an efficient computation.  The multiscale
method provides a well-adapted spatial representation of the solution, and the
dual solution will detect time-domains where the solution is stationary. In
these domains, the equidistribution strategy will choose large timesteps.

{\bf Target Functional.}
Now we set up the target functional. The functional $J(U)$ is chosen as a
weighted average of the normal force component exerted on the bump and at the
boundaries before and behind the bump:
\begin{align}\label{eq:func_p}
  J(U) =\sum_{i=1}^7 \int_0 ^T \int_{\kappa_{i}} p \,\psi_{i} (x,y) ds
\end{align}
with
\begin{align*}
  \kappa_{i} &= \left\{
  (x,y)\in \Gamma \,:\, x\in [x_i-0.25, x_i+0.25]\right\}\\
  \psi_{i}(x) &= (x-(x_i-0.25))^2(x+(x_i+0.25))^2/0.25^4,\quad x\in \kappa_i. 
\end{align*}

Here $\Gamma$ is only the bottom part of $\Gamma$. In all computations presented
in this section the functional \eqref{eq:func_p} is chosen, which is the
pressure averaged at several points at the bump in front and behind the bump.
The $x$-coordinates of these points at the bottom are $x_i$ = -3, -2, -1, 0, 1,
2, 3. At each of these points $x_i$ a smooth function $\psi_{i}$ is given with
support $x_i-0.25, x_i+0.25$. This functional measures the pressure locally.

In Figure \ref{channel-stat3} we show a time-sequence of the instationary test 
case computed with uniform $\cfl=1$ on an adaptive grid with finest level $L=5$. 
Note that the perturbation entering at the left boundary and moving through the 
boundary is resolved very well.

\section{Computational results}
\label{section.numexp_fully_implicit}

\subsection{Numerical strategies}

We gave an outline of the adaptive method in Section \ref{section.numexp}. Now we will present numerical computations, where we compare some variations of the adaptive concept. 

  The first strategy is the one we proposed in \cite{SteinerNoelle}:
  We first compute a forward solution on a coarse grid ($L=2$) and solve 
  the adjoint problem on the adaptive grid of the forward solution.
  Then we use the information of the error representation based on the dual solution to determine a new sequence of timesteps. This sequence is used in the computation of the forward solution on a grid with finest level $L=5$, where we additionally restrict the $\cfl$ number from below. 
  We compare the results of the time adaptive strategy with uniform timestep distributions.

In the second strategy we modify the fully implicit timestepping strategy and introduce a mixed {\em
implicit/explicit} approach. The reason is that implicit
timesteps with $\cfl<5$ are not efficient, since we have to solve
a nonlinear system of equations at each timestep. Thus, for $\cfl<5$,
the new implicit/explicit strategy switches to the cheaper and less dissipative
explicit method with $\cfl=0.5$.

We want to compare these strategies with respect to the following main aspects:
\begin{itemize}
\item What is the quality and what are the costs determining the adaptive
  timestep sequence from computations on the coarse grid?
\item Is the predicted adaptive timestep sequence well-adapted to the solution
  on the fine grid?
\item How is the solution affected if we use uniform timesteps larger than
  the predicted adaptive timestep sequence?
\end{itemize}

In order to quantify the results we have to compare with a reference solution.
Since the exact solution is not available we perform a computation with $L=5$
refinement levels using implicit timestepping with $\cfl= 1$ amd explicit timestepping with $\cfl= 0.5$. 
\subsection{Strategy I: Fully implicit computational results}\label{sec.numadjoint}
\subsubsection {Adjoint indicator and adaptive timesteps}

Now we use the error representation on finest level $L=2$ for a new time
adaptive computation on finest level $L=5$. 

\begin{figure}[ht]  
  \begin{center} {\vspace{-25mm}
  {\includegraphics[width=0.9\textwidth]{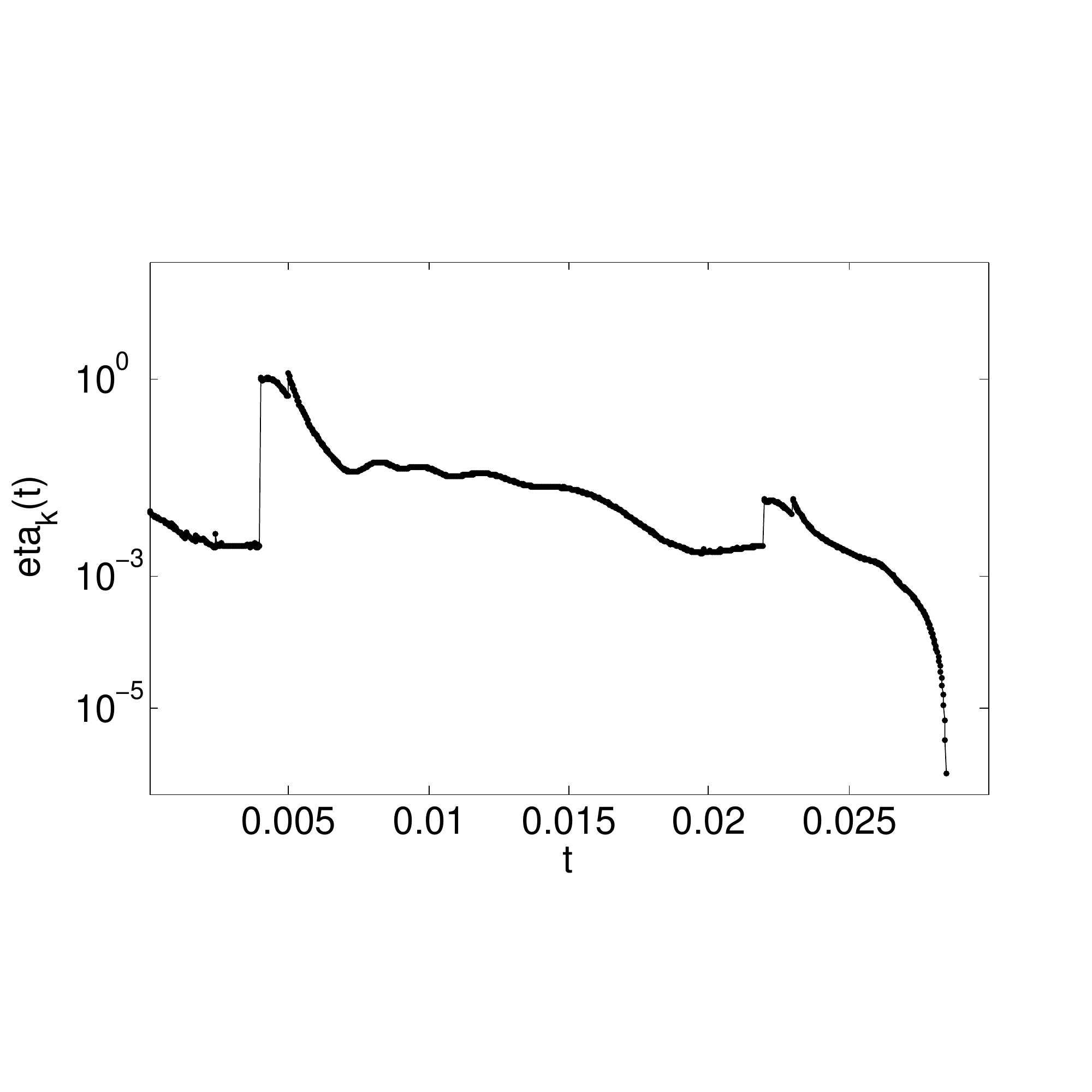}}\vspace{-4cm}
  {\includegraphics[width=0.9\textwidth]{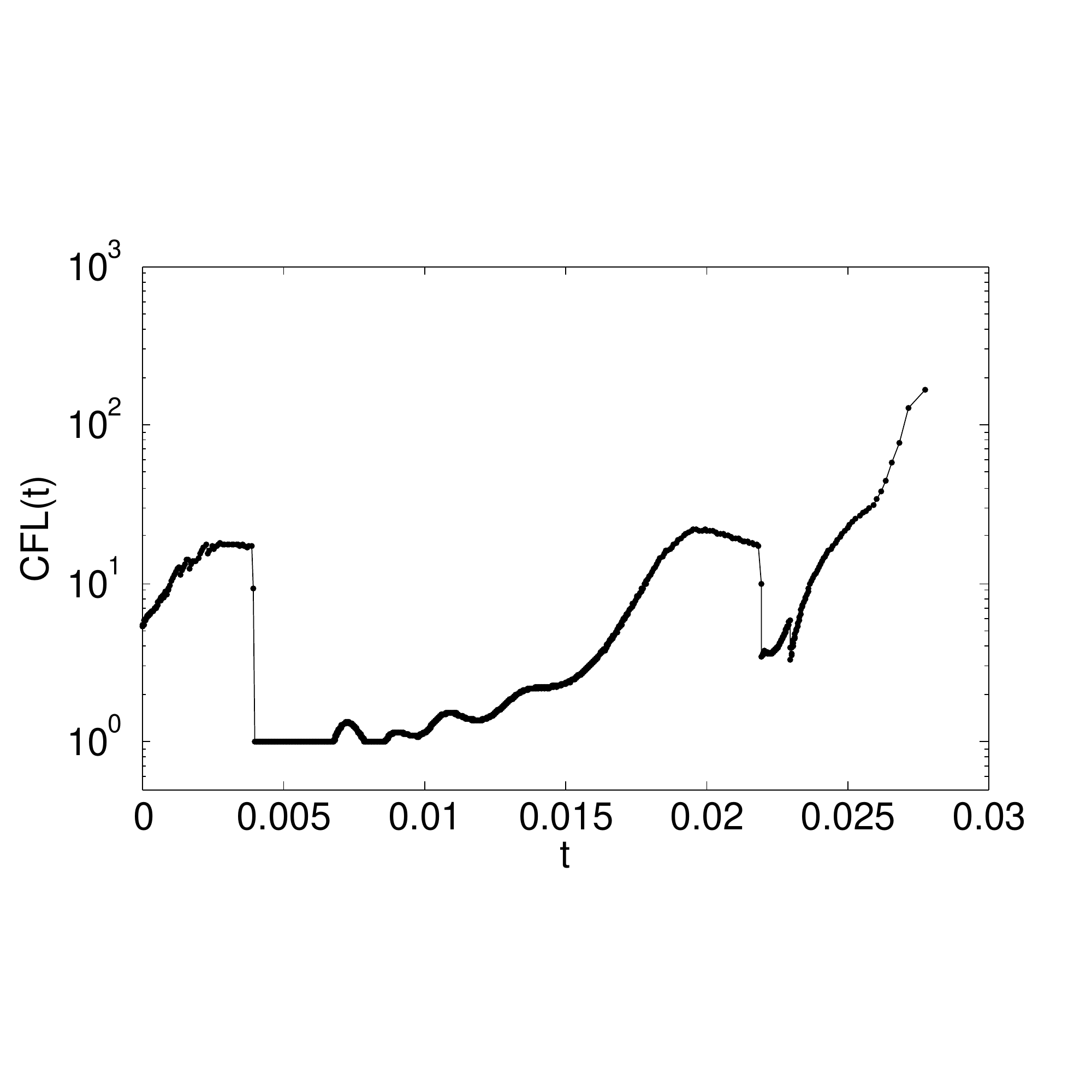}}\vspace{-2.5cm}}
  \end{center}
  \caption{Time component of the error representation $\ebkn$(top) and new
    timesteps with $\cfl$ restriction from below $\cfl(t_n)$ (bottom).
  \label{fig.vgl_euler2d_1}}
\end{figure}

The first computation is done on a mesh with finest level $L=2$. We compute
until time $T= 0.0285 s$, which takes 1000 timesteps with $\cfl=1$. We use the
results of the error representation of this computation to compute a new
timestep distribution. The forward problem takes $329 s$ and the dual problem
including the evaluation of the error representation $619 s$ on an Opteron 8220
processor at 2.86 GHz. The total computational costs are $948 s$, and in memory
we have to save 1000 solutions (each timestep) of the forward problem which
corresponds to 48 MB (total). This gives us a new sequence of adaptive timesteps
for the computation of level $L$ = 5. The error indicator and the new timesteps
are presented in Figure \ref{fig.vgl_euler2d_1}. In time intervals where the
solution is stationary, i.e. at the beginning, and after the perturbations have
left the computational domain, the timesteps are large. In time intervals where
the solution is instationary we get well-adapted small timesteps.

Then we use the adaptive timestep sequence for a computation on level $L=5$ and
compare it with a uniform in time computation using $\cfl= 1$. The uniform
computation needs 8000 timesteps and the computational time is $21070 s$. The
time adaptive solution is computed with 2379 timesteps and this computation
takes 9142$s$.

In Figure \ref{fig.vgl_euler2d_4} we show a sequence of plots of the uniform, $\cfl=1$ computation on an adaptive spatial grid with finest level $L=5$. In Figure \ref{fig.vgl_euler2d_4} we compare the pressure distribution at the bottom boundary of the uniform solution and the time adaptive solution at several times. The two solutions on level $L=5$ match very well.

\begin{remark}
In \cite{SteinerMuellerNoelle} we also did some comparisions of the adjoint indicator with some ad hoc indicators, which estimate the variation of the solution from one timestep to the following. While these indicators also detect whether the solution is stationary or not, they lead to computations which are in general more expensive but not more accurate. 
In particular, most timesteps
are smaller than in the case with adaptation via adjoint problems.

 One
advantage may be that we do not need to compute a dual solution, which makes the
computation of the variation indicator less expensive. But this is only a small
advantage, since we compute the error indicators on a coarse mesh.
On the other hand the ad-hoc indicator is more straight forward to implement.

Another approach was to choose the maximum jump
of the solution in one cell, both weighted and not weighted with the size of the
cell. This was an approximation to the $L^\infty$-norm. This indicator is not
very useful, since it turned out to be highly oscillatory. 
\end{remark}

\subsubsection{Comparison with Uniform timesteps}\label{sec.numuniform}

In many instationary computations where no a-priori information is known, one
reasonable choice is to use uniform $\cfl$ numbers. Therefore we compared
the computation with adaptive implicit timesteps with implicit computations
using uniform $\cfl$ numbers. In Section~\ref{sec.numadjoint} we have already
done a computation with uniform $\cfl$ number, $\cfl=1$, on a grid with $L=2$,
to get timestep sizes for an adaptive computation on a grid with $L=5$.  As a
reference solution we also computed with uniform $\cfl$ number, $\cfl =1$, a
solution of the problem on a grid with $L=5$. Now we compare these computations
with computations using higher uniform $\cfl$ numbers.

\begin{figure}[t]
\begin{center}\vspace{-2cm}\hspace{-1.3cm}
\includegraphics[width=0.6\textwidth]{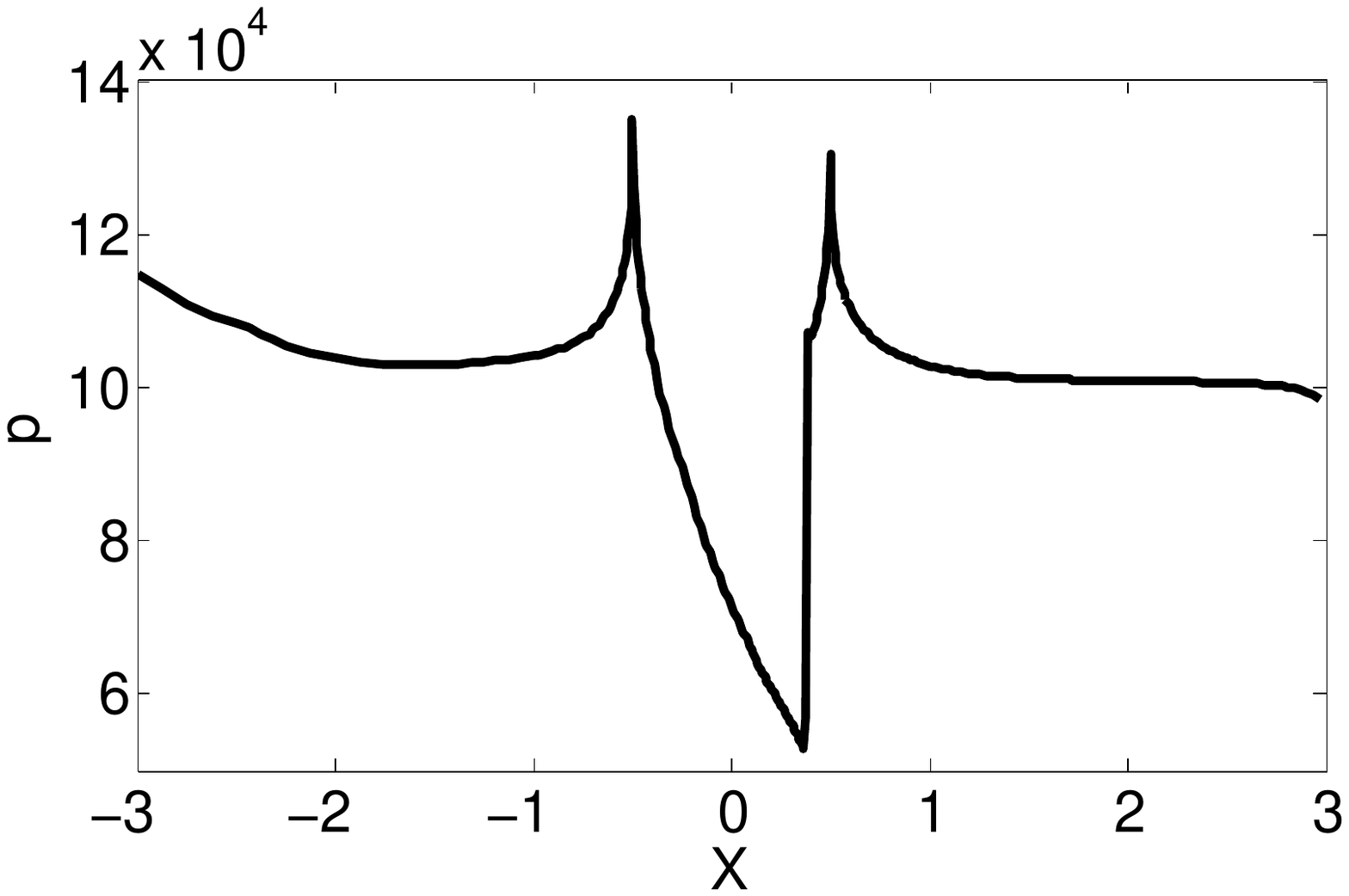}\hspace{-1.3cm}
\includegraphics[width=0.6\textwidth]{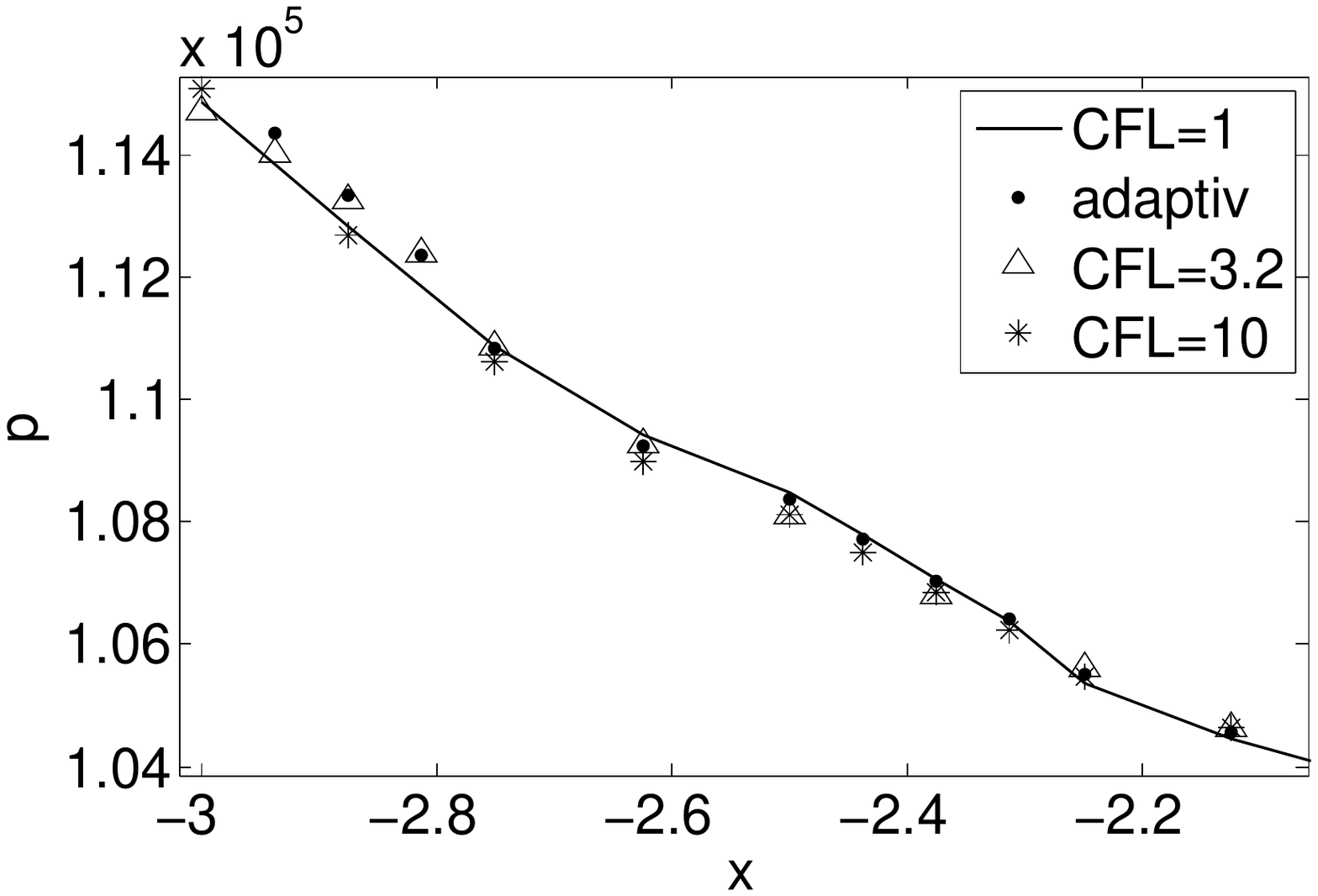}\\
\vspace{-5.5cm}\hspace{-1.3cm}
\includegraphics[width=0.6\textwidth]{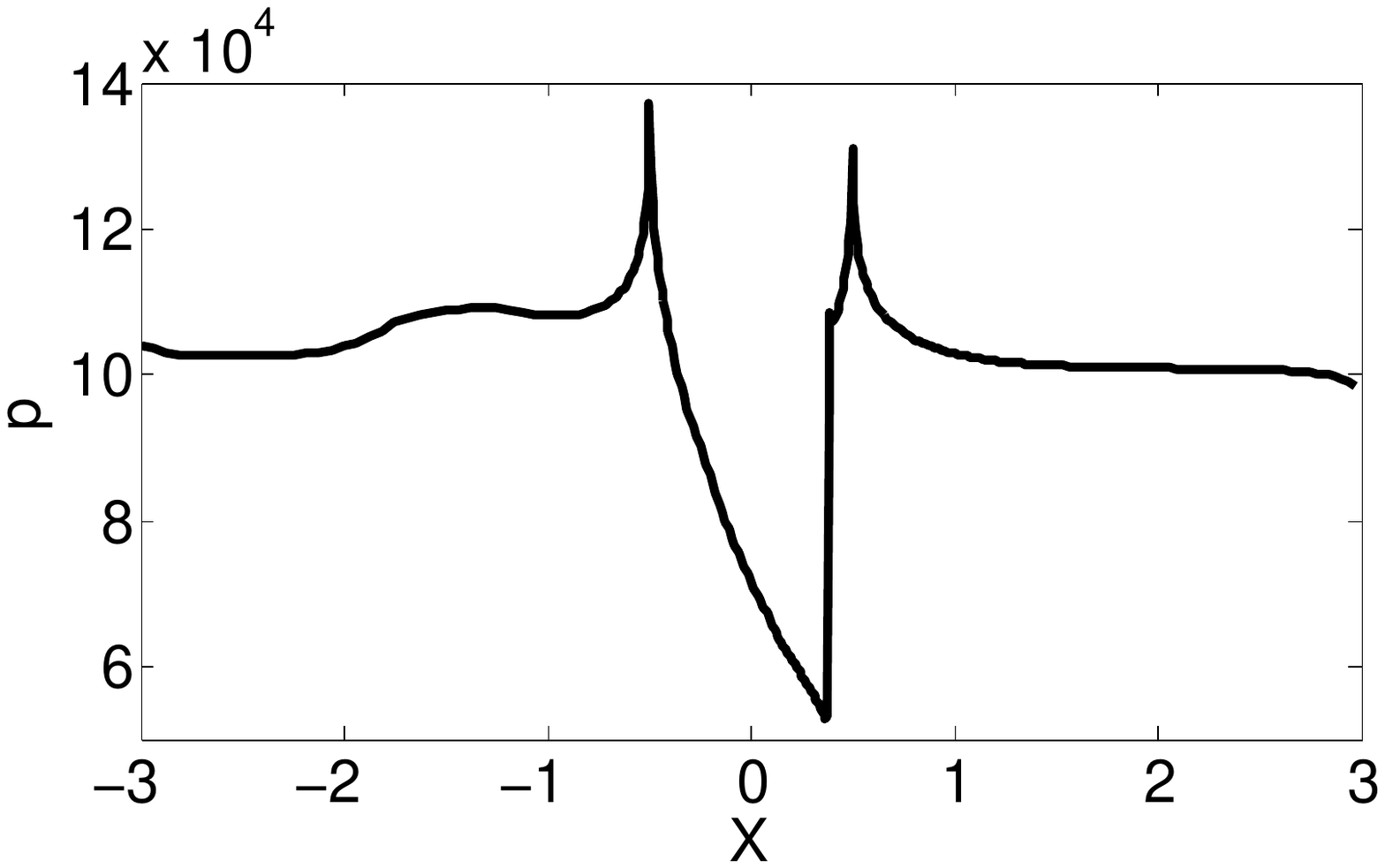}\hspace{-1.3cm}
\includegraphics[width=0.6\textwidth]{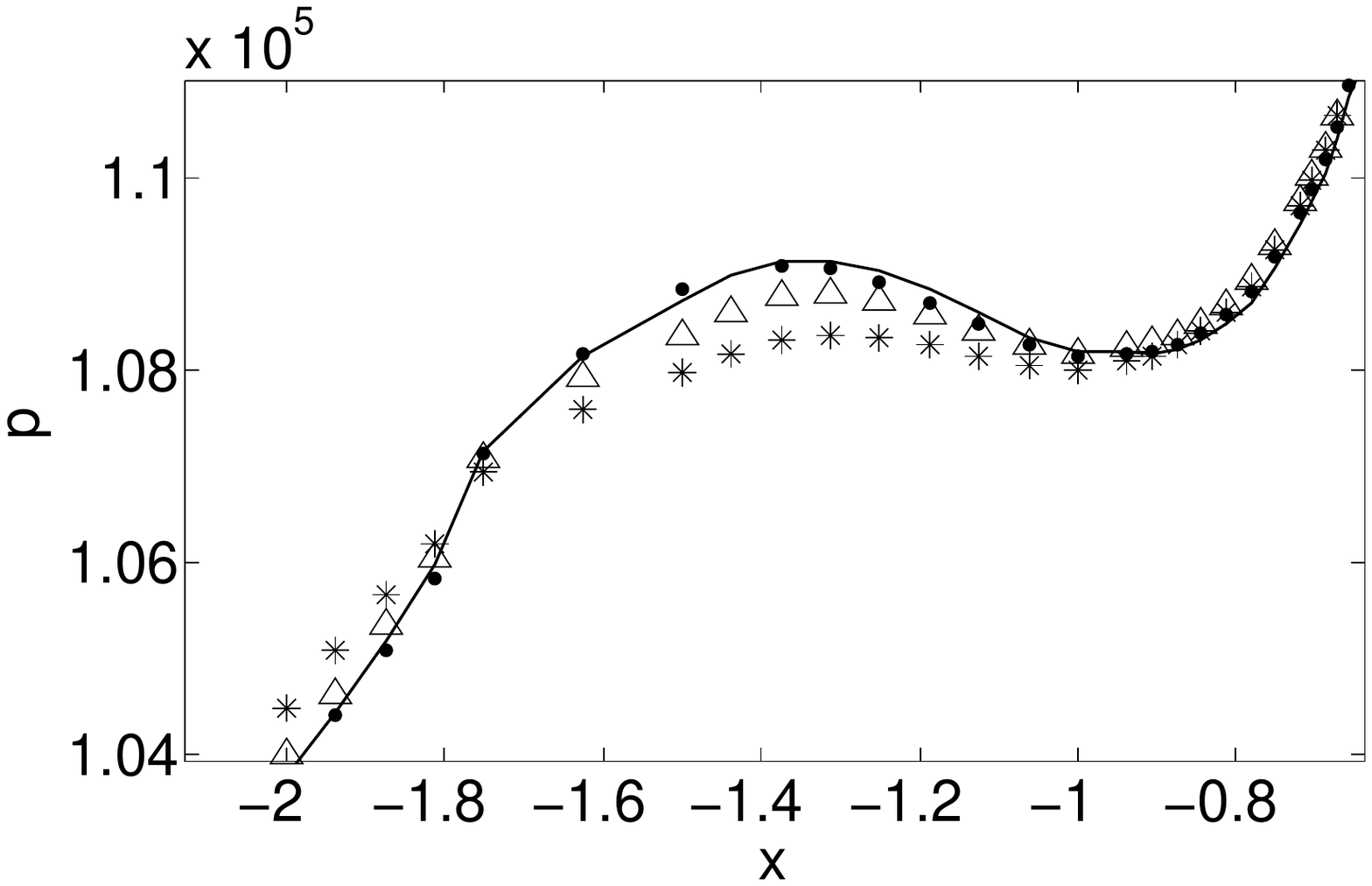}\hspace{-1.3cm}\\
\vspace{-5.5cm}\hspace{-1.3cm}
\includegraphics[width=0.6\textwidth]{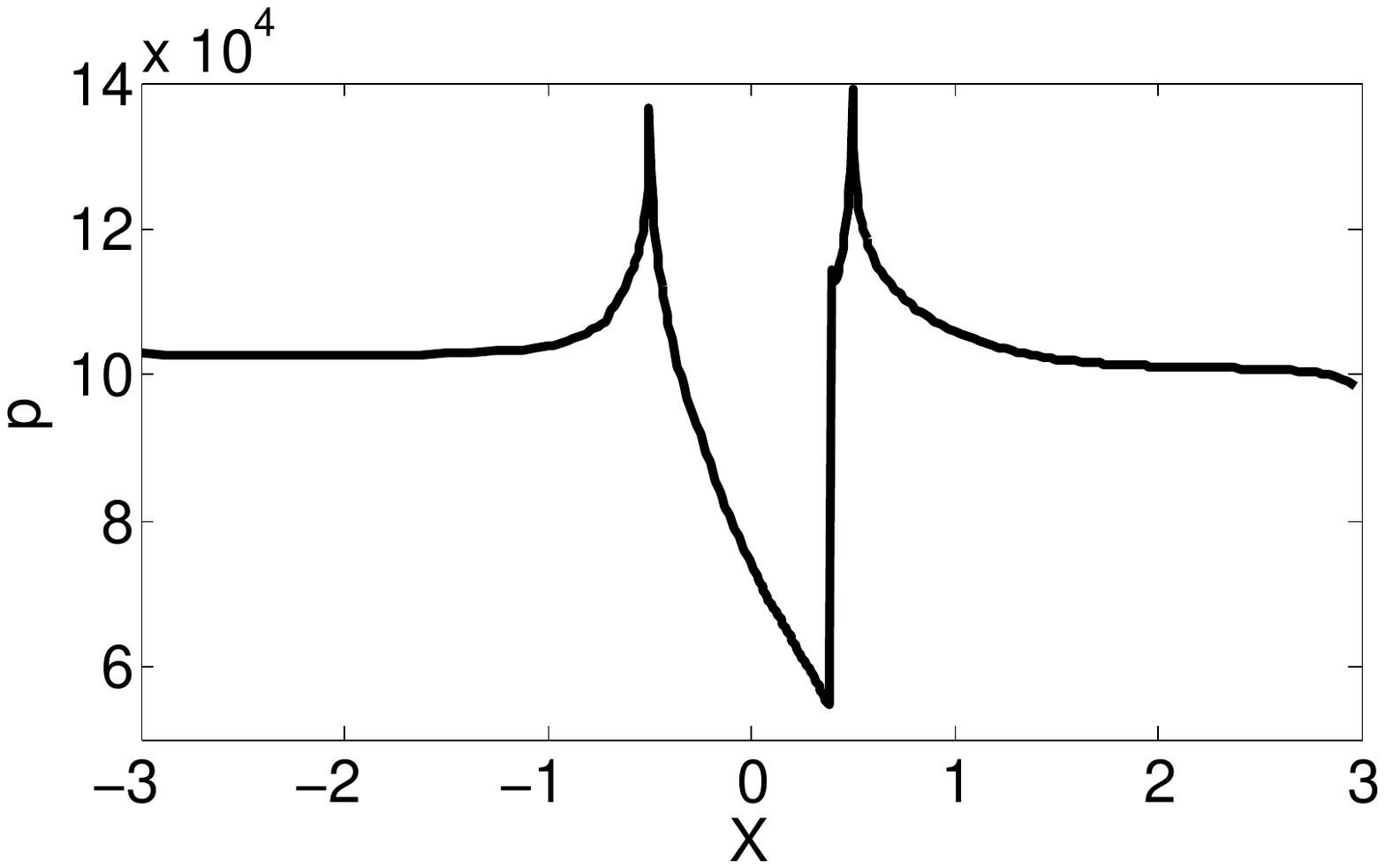}\hspace{-1.3cm}
\includegraphics[width=0.6\textwidth]{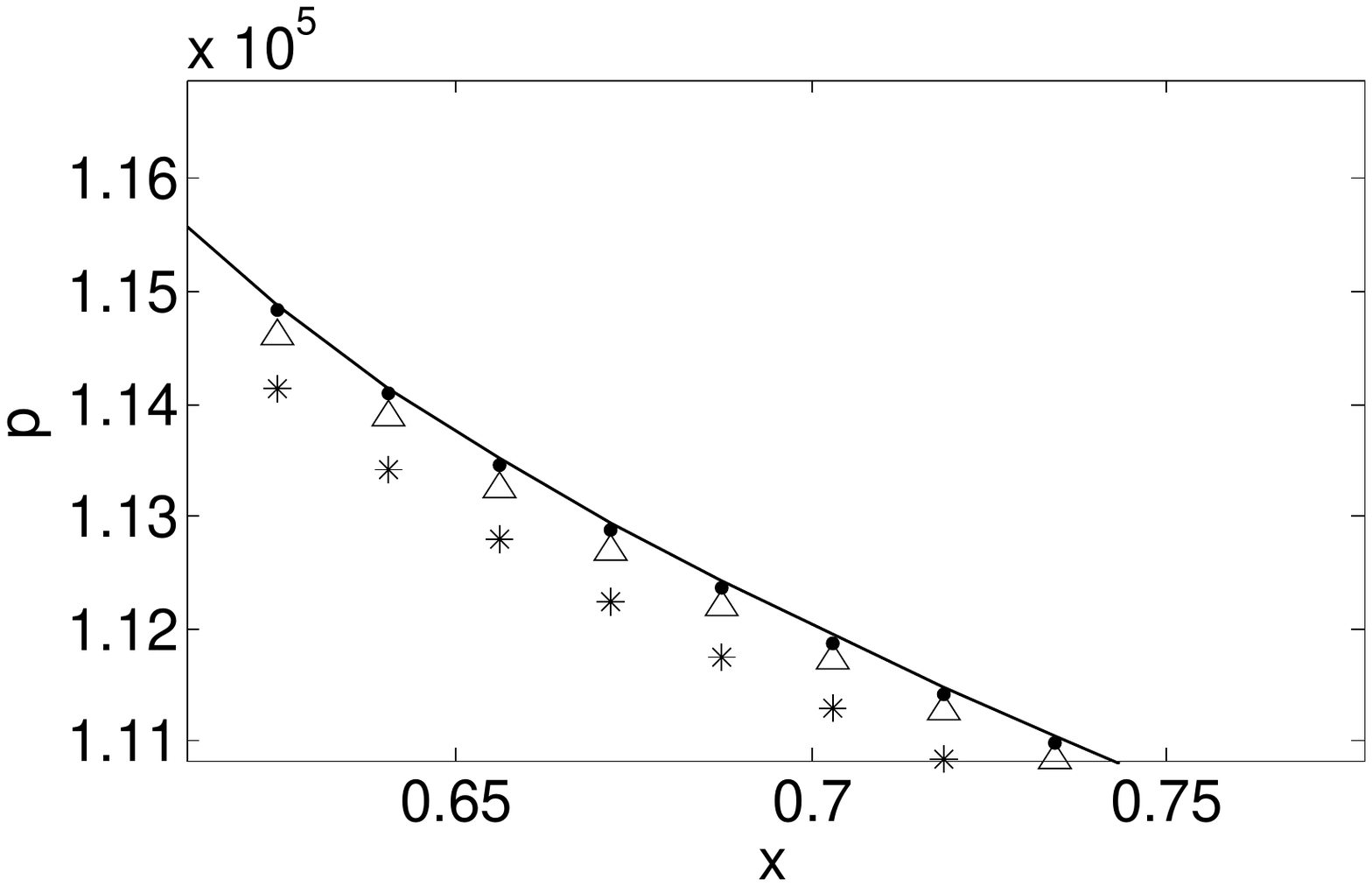}\hspace{-1.3cm}\\
\vspace{-5.5cm}\hspace{-1.3cm}
\includegraphics[width=0.6\textwidth]{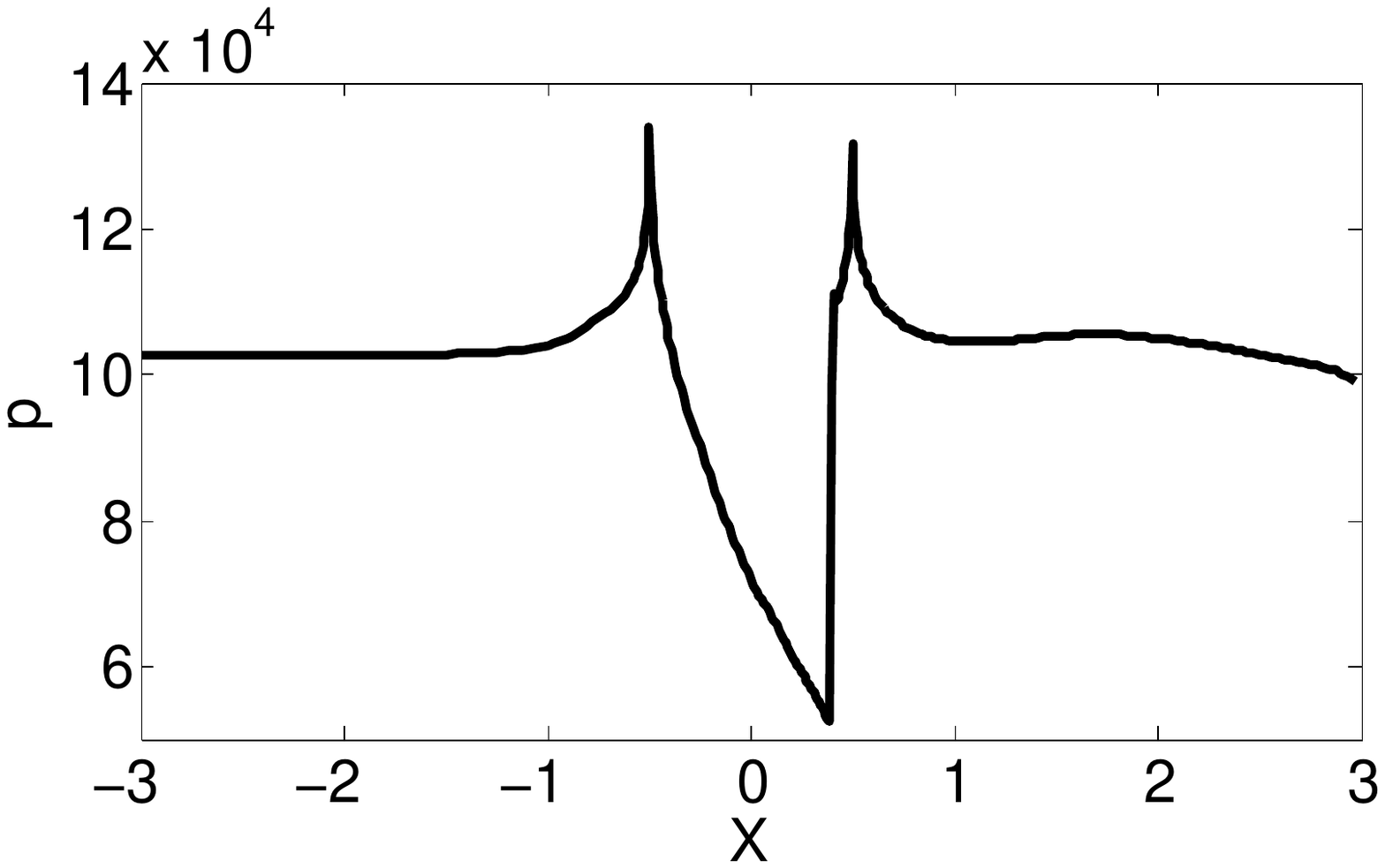}\hspace{-1.3cm}
\includegraphics[width=0.6\textwidth]{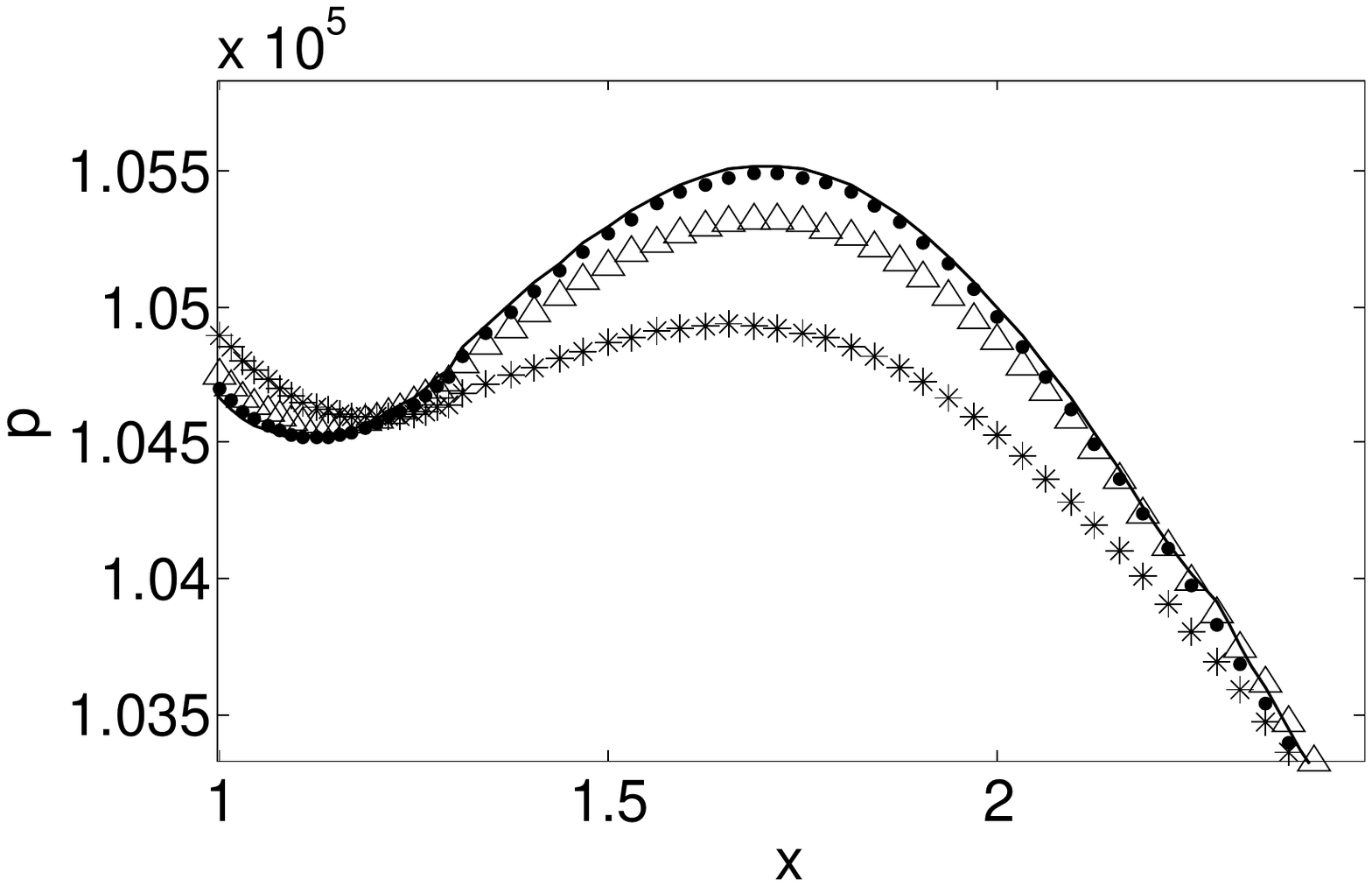}\hspace{-1.3cm}\\
\vspace{-3cm}
\caption{Reference solution for 2D Euler equations: Left: adaptive spatial grid (finest level $L=5$), uniform time steps ($\cfl=1$). Pressure $p$ at bottom boundary at times $t$=0.005002, 0.007125, 0.009990, 0.011975, from top to bottom. Right: Zoom tracing the perturbations. Comparison of uniform timesteps with $\cfl = 1,\;3.2,\;10$ and time-adaptive strategy.} \label{fig.vgl_euler2d_4}
\end{center}
\end{figure}
First we choose a uniform $\cfl$ number of approximately 3.2, which corresponds
to 2500 timesteps. This equals roughly the number of timesteps in the
adaptive method, and hence it should give a fair comparison. The uniform
computation takes about $11509 s$, more than the $9142 s$ of the adaptive
computation (see Table~\ref{table.comp2}).  In the uniform computation most of
the timesteps are more expensive, since they need more Newton steps, and more
steps for solving the linear problems. This shows that the understanding of the
dynamics of the solution pays directly in the nonlinear and linear solvers.
Moreover, it can be seen from Figure~\ref{fig.vgl_euler2d_4} that the quality of
the solution is considerably worse than for the adaptive computation.

Another computation with $\cfl$ number 10 takes only $3290 s$. However, as can
be seen from Figure~\ref{fig.vgl_euler2d_4} the solution is badly approximated:
In the beginning of the computation the solutions of the different methods match
very well, which means that the inflow at the boundary is well-resolved. As time
goes on, the solutions differ more and more. After the perturbation has passed
the bump, the perturbations differ strongly. Only the time-adaptive method
approximates the reference solution ($\cfl=1$) closely.

\begin{remark}
Table~\ref{table.comp2} gives an overview of the CPU time and the number of
Newton iterations and linear iterations for the computations. The costs for the
computation of the indicator on level
$L=2$ are very low compared to the costs of computations on level $L=5$.  An
adaptive computation including the computation of the indicator, i.e. $329 s +
619 s + 9142 s$, is cheaper than the computation using uniform $\cfl$ number,
e.g. $\cfl=1$, that needs 21070$s$. Even the computation with $\cfl=3.2$ is more
expensive than the time-adaptive computation, but leads to worse results, see
Figure \ref{fig.vgl_euler2d_4}.

Table~\ref{table.comp2} shows that the CPU time is roughly proportional to the
number of Newton iterations and not to the number of timesteps or linear solver
steps. The CPU time is about $2.5 s$ per Newton iteration. This means that we
have to minimize the number of Newton iterations in total to accelerate the
computation. This is done very efficiently by the time-adaptive approach. For a
large range of $\cfl$ numbers from 1 to more than 100, it needs only one or two
Newton iterations per timestep, without sacrificing the accuracy, see \cite{SteinerMuellerNoelle}.

In this example, the stationary time regions are not very large compared to the
overall computation. If the stationary regions were larger, the advantage of the
time adaptive scheme would be even more significant.

\begin{table}[ht]  
\begin{center}
\begin{tabular}{l|c|c|c|c}
&CPU [s]& timesteps & Newton steps & linear steps \\ \hline
adaptive timesteps & 9142& 2379 &3303&16875\\
uniform $\cfl=1$ & 21070&8000&8282&32630\\
uniform $\cfl=3.2$ & 11509&2500 &4904&26895\\
uniform $\cfl=10$ & 3290&800&1600& 13500\\
\end{tabular}\\[1ex] 
\caption{Performance for computations on $L=5$ using different fully implicit
timestepping strategies.}\label{table.comp2}
\end{center}
\end{table}

\end{remark}

\subsection{Strategy II: Explicit-implicit computational results}
\label{section.numexp_explicit_implicit}

\begin{figure}
  \begin{center}
  \includegraphics[width=0.75\textwidth]{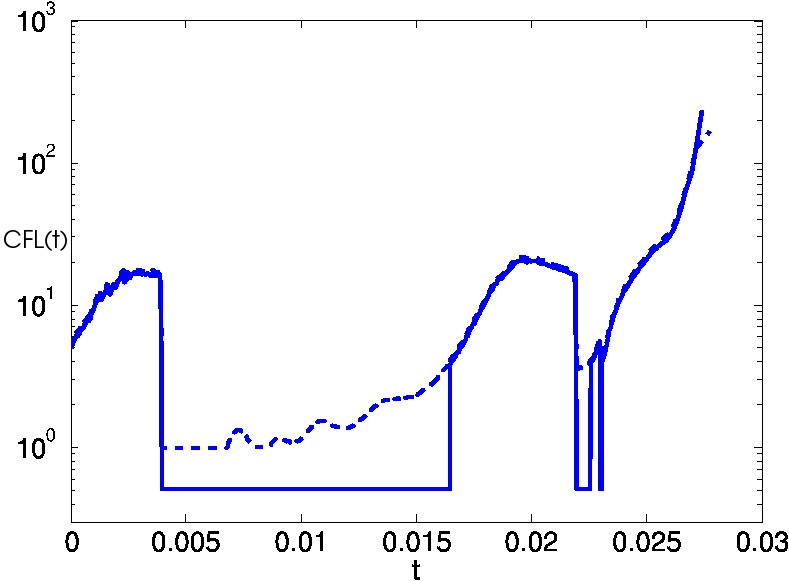}
  \caption{2D Euler equations, comparison of timestep sequence derived from error
  representation for implicit computation (dashed line) and for mixed
  explicit-implicit computation  (bold line).}
  \label{fig.vgl_ex-im_CFL}
  \end{center}
\end{figure}

Now we modify the fully implicit timestepping strategy and introduce a mixed {\em
implicit/explicit} approach. The reason is that implicit
timesteps with $\cfl<5$ are not efficient, since we have to solve
a nonlinear system of equations at each timestep. Thus, for $\cfl<5$,
the new implicit/explicit strategy switches to the cheaper and less dissipative
explicit method with $\cfl=0.5$. The timestep sequence is shown in 
Figure~\ref{fig.vgl_ex-im_CFL}. Of course, we could choose variants of
the thresholds $\cfl=0.5$ and 5.

\begin{figure}
  \begin{center}
  \includegraphics[width=0.85\textwidth]{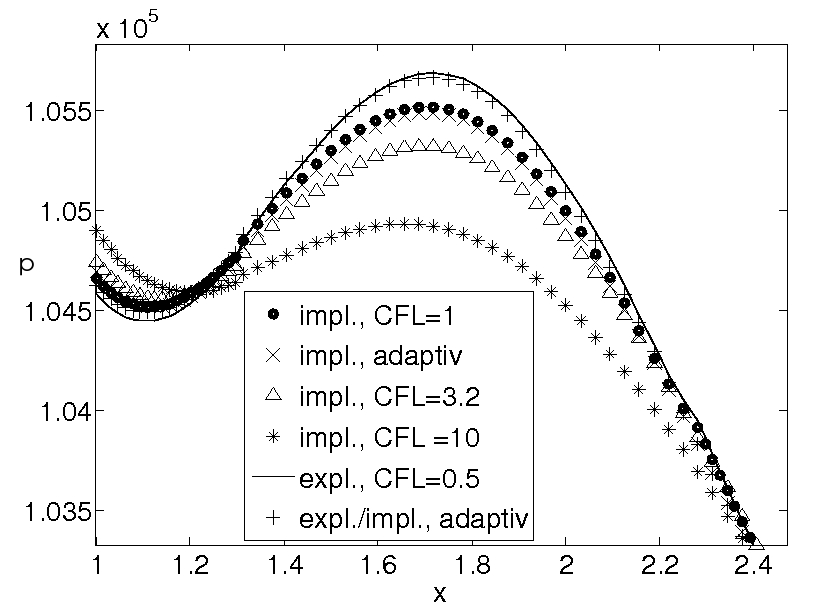}
  \caption{2D Euler equations, comparison of solutions in time on adaptive grid
  with finest level $L=5$, pressure $p$ at the bottom boundary, zoom of the
  perturbation at time $t$=0.011975}
  \label{fig.vgl_ex-im}
  \end{center}
\end{figure}

As we can see in Table~\ref{table.comp3}, the new strategy requires 5802
timesteps, where $95\%$ are explicit. The CPU time of $7730 s$ easily beats
the fully explicit solver ($18702 s$), and is also superior to the fully
implicit adaptive scheme ($ 9142 s$, see Table \ref{table.comp2}). 
The computational results are presented in Figure \ref{fig.vgl_ex-im}. The
results of the combined explicit-implicit strategy are very close to the results
of the fully explicit method, and far superior to all fully implicit methods. 
Note that the explicit scheme serves as reference solution, since it is
well-known that it gives the most accurate solution for an instationary problem.

\begin{table}[ht]  
\begin{center}
\begin{tabular}{l|c|c|c}
&CPU [s]& timesteps & timesteps \\
&       & (total)   & (implicit)\\ \hline
explicit $\cfl=0.5$ & 18702&16000&-\\
adaptive expl.-impl. & 7730&5802&259\\
\end{tabular}\\[1ex] 
\caption{Performance for computations on $L=5$ using fully explicit timesteps
and the time-adaptive explicit-implicit strategy.}
\end{center}
\label{table.comp3}
\end{table}

\section{Conclusion}
\label{section:conclusion}

In the series of works \cite{Steiner_2008_diss, SteinerNoelle, SteinerMuellerNoelle} reviewed in this paper, explicit and 
implicit finite volume solvers  on adaptively refined
quadtree meshes have been coupled with adjoint techniques to control the
timestep sizes for the solution of weakly instationary compressible inviscid
flow problems. 

For the 2D Euler equations we study a test case for which the
time-adaptive method  does reach its goals: it separates stationary regions and
perturbations cleanly and chooses just the right timestep for each of them. The
adaptive method leads to considerable savings in CPU time and memory while
reproducing the reference solution almost perfectly.

Our approach is based upon several analytical ingredients:
In Theorem~\ref{theorem:error_representation_h} we state a complete error
representation for nonlinear initial-boundary-value problems with characteristic
boundary conditions for hyperbolic systems of conservation laws, which includes
boundary and linearization errors. Besides building upon well-established
adjoint techniques, we also add a new ingredient which simplifies the
computation of the dual problem \cite{SteinerNoelle}. We show that it is
sufficient to compute the  spatial gradient of the dual solution, $w=\nabla
\vp$, instead of the dual solution $\vp$ itself. This gradient satisfies a
conservation law instead of a transport equation, and it can therefore be
computed with the same algorithm as the forward problem. For discontinuous
transport coefficients, the new conservative algorithm for $w$ is more robust
than transport schemes for $\vp$, see \cite{SteinerNoelle}. Here we also derive
characteristic boundary conditions for the conservative dual problem, which we
use in the numerical examples in Section~\ref{section.numexp_fully_implicit}.

In order to compute the adjoint error representation one needs to compute a
forward and a dual problem and to assemble the space-time scalar product
\eqref{eq:localized_error_indicator}. Together, this costs about three times as
much as the computation of a single forward problem. In our application, the
error representation is computed on a coarse mesh ($L=2$), and therefore it
presents only a minor computational overhead compared with the fine grid
solution ($L=5$). In other applications, the amount of additional storage and
CPU time may become significant.

For instationary perturbations of a 2D stationary flow over a bump, we have 
implemented and tested both a fully implicit and a mixed
explicit-implicit timestepping strategy. The explicit-implicit approach switches
to an explicit timestep with $\cfl=0.5$ in case the adaptive strategy suggests
an implicit timestep with $\cfl<5$. Clearly, the mixed explicit-implicit
strategy is the most accurate and efficient, beating the adaptive fully implicit
in accuracy and efficiency, the implicit approach with fixed $\cfl$ numbers in
accuracy, and the fully explicit approach in efficiency.

\bibliographystyle{styles/spmpsci}
\bibliography{christina_lit2}
\end{document}